\documentclass[11pt,a4paper,notitlepage]{article}
\usepackage{latexsym, amsmath, verbatim, graphics, epsfig, amssymb, color}
\newtheorem{theorem}{Theorem}[section]
\newtheorem{lemma}[theorem]{Lemma}

\newtheorem{definition}[theorem]{Definition}
\newtheorem{remark}[theorem]{Remark}

\begin{document}
\begin{center}
\textbf{Exponential Stability of Higher-Order Fractional Neutral Stochastic Differential Equation via Integral Contractors}
\vskip 0.3cm   
Dimplekumar Chalishajar${^{}}$\footnote{Corresponding  author.
	Email: chalishajardn@vmi.edu (Dimplekumar Chalishajar). Tel: +1-540 817 8105.}, K. Dhanalakshmi${{}^{2}}$,
K. Ramkumar${{}^{3}}$,
K. Ravikumar${{}^{4}}$
${}^{1}$Department of Applied Mathematics, Virginia Military Institute(VMI), Lexington,\\ VA 24450, USA.\\
\vskip 0.1cm ${}^{2}$Department of Mathematics,
The Gandhigram Rural Institute (Deemed to be University),\\
Gandhigram - 624 302, Tamil Nadu, India.\\
${}^{3,4}$Department of Mathematics, PSG College of Arts and Science, Coimbatore, 641 014, India.\\  ${{}^{1}}$chalishajardn@vmi.edu,${{}^{2}}$dlpriya20@gmail.com, ${{}^{3}}$ramkumarkpsg@gmail.com, ${{}^{4}}$ravikumarkpsg@gmail.com
\end{center}
\textbf{Abstract:} The existence, uniqueness and exponential stability results for mild solutions to the fractional neutral stochastic differential system are presented in this article. To demonstrate the results, the concept of bounded integral contractors is combined with the stochastic result and sequencing technique. In contrast to previous publications, we do not need to specify the induced inverse of the controllability operator to prove the stability results, and the relevant nonlinear function does not have to meet the Lipschitz condition. Furthermore, exponential stability results for neutral stochastic differential system with Poisson jump have been established. Finally, an application to demonstrate the acquired results is discussed. This paper extends work of Chalishajar et al. \cite{r4} and Renu Chaudhary et al. \cite{r3}. 
\vspace{0.2cm} \vspace{0.2cm}

\section{Introduction}
The mathematical community began paying more and more attention to fractional calculus, or FC. The majority of early work was on creating analytical formulations to address particular mathematical issues. The expansion of definitions for fractional operators, such as the integral representation (Liouville, Riemann, and Hadamard) and the convergent series representation (Grunwald and Letnikov), was the most direct outcome of the fast increasing interest in fractional continuity (see the monograph \cite{r9, r11, r12, r14, r16}). The majority of early research remained focused on the creation of the mathematical framework and the integration of these operators into ordinary and partial differential equations, despite the fact that these early studies had highlighted the fascinating role that FC can play when modeling complex processes in physical systems. It was not until the latter part of the 20th century that the idea of FC began to spread beyond of mathematics. One application area that has grown quite quickly is the simulation of intricate physical events. In fact, viscoelastic effects, nonlocal behavior, anomalous and hybrid transport, fractal media, and even control theory were among the numerous early physical modeling applications of FC \cite{r6, r10, r20}.

\noindent
The integer-order situation has been the focus of most of these models' investigations. When attempting to use a typical integer-order explanation to describe the dynamics of a particle event, researchers resort to fractional calculus. Fractional calculus has drawn a lot of interest from specialists in science and engineering as it may be utilized in a wide variety of mathematical models in these domains. This is due to the fact that fractional calculus allows for the exploration of the fractional operators' utility in a variety of contexts. Three fractional operators—Caputo, Caputo-Fabrizio, and Atangana-Baleanu—are often used in recent scholarly research. In recent research, a fractional model is offered for the banking data using the framework developed by Caputo fractional derivative, with model parameters generated by least-squares curve fitting and a larger range of Caputo-Fabrizio operators. We elucidate the several acronyms that will be employed in this review to denote the various categories of fractional-order operators \cite{r9, r11, r12}, including:
\begin{enumerate}
\item "CO" operators stand for single constant-order operators.
\item "DO" stands for distributed-order operators (with constant order distribution).
\item   "VO" operators stand for variable-order operators. Although VO operators might be distributed or single, when we use the abbreviation "VO," we mean single variable-order operators alone. 
\item "DVO" operators stand for distributed-variable-order operators, which will be covered later.
\end{enumerate}
The groundwork for distributed-order fractional calculus (DOFC) was laid by Caputo's groundbreaking research on dissipative elastodynamics. In these investigations, a parallel series of fractional-order derivatives was used to generalize the viscoelastic stress–strain constitutive laws. Originally, this operator was known as the "mean fractional-order derivative." The many viscoelastic models that may be retrieved from the multi-term DO law in the following examples further highlight the efficacy of the DO approach:

\noindent
\textbf{Kelvin-Voigt Models:}  The DO analogue of the Kelvin–Voigt model is obtained for the choice of $\phi_{\sigma} =  \delta(\gamma)$ and  $\phi_{\epsilon}  = \kappa^{\gamma}$.

\noindent
\textbf{Maxwell Models:} The fractional-order Maxwell model of viscoelasticity can be obtained for $\phi_{\sigma}= \delta (\gamma)+ \kappa^{\alpha}\delta(\gamma-\alpha)$  and  $\phi_{\epsilon}= E_{\infty} \kappa^{\beta} \delta(\gamma-\beta)$. Note that, assuming  $\alpha = \beta  $ in the fractional Maxwell model, allows recovering the fractional Zener model.

\noindent
\textbf{Zener Model:} Wave motion in fractional To create viscoelastic media of the Zener type, select   $\phi_{\sigma}=  \phi_{\epsilon}= \delta (\gamma)+ \kappa^{\alpha}\delta(\gamma-\alpha)$. Similarly, the choice of $\phi_{\sigma}= \delta(\gamma)+ \frac{a}{b}  \delta(\gamma- (\alpha-\beta))$ and $\phi_{\epsilon} = a \delta(\gamma-\alpha) + c\delta(\gamma-\eta)+\Big(\frac{ac}{b}  \Big) \delta (\gamma-\alpha-\eta+\beta)$, also produces a fractional version of the spring-and dashpot-based conventional Zener model.\\ 
\noindent
In actuality, these stochastic dynamic systems depend on both the past and current states. For example, \cite{r5, r6, r15} describes these systems in terms of stochastic functional differential equations. Stability analysis, as a popular issue in stochastic dynamical systems research, has caused a significant deal of anxiety (see monographs \cite{r5} and \cite{xm}). There is now a wealth of material about SFDE stability. One of the fundamental characteristics of the world is stochasticity, and in practical system application, stability is of utmost importance. As everyone knows, stability is crucial for a system with real-world application experience. Numerous theories of system stability have been proposed by scientists on real-world requirements.  For stochastic systems, stability analysis is crucial (see \cite{r6, r7, r15}).\\
\noindent
For abstract differential equations, the existence results are established by employing the successive approximation method. Also, the semigroup theory of bounded linear operators is firmly identified
to find solvability of abstract differential equations within infinite-dimensional spaces \cite{r13}. By using an infinitesimal generator $ \mathcal{A} $ of $ C_{0}-$semigroup of bounded linear operators $ C_{\alpha} (\iota) $ on Hilbert spaces, one can construct a mild solution for the given system.\\
\noindent
The concept of contractors by Altman (1978) \cite{r1} has been extended in the case of random operators on Banach spaces. The nonlinear term is assumed to have an integral contractor which is a weaker condition than the Lipschitz’s continuity. The existence of the system is guaranteed by the integral contractor but not the uniqueness \cite{r3}. We rely
on the definition of regularity to check the uniqueness of the system’s mild solution. This technique can be implemented to all dynamical, stochastic, and fractional order systems, for more details on integral contractors, readers refer the articles \cite{r2, r3, r8, r11, r20}.\\

\noindent {\bf{Novelty:}}
\begin{enumerate}
	\item Stability of the stochastic system using integral contractor and regularity, has been investigated. No research has been reported based on this new technique. 
	\item Sufficient conditions for the higher order fractional stochastic system is studied  when the nonlinear term is not Lipschitz continuous.	
	\item DOFC used in the system has a multiple applications in Kelvin-Voigt models, Maxwell models, Zener model, etc. System considered here is the extended form of the these models. 
\end{enumerate}

\noindent The article framework may be divided into the following categories: In Sect. 2, we give some definitions of fractional calculus, semigroup theory, some associated notations and essential lemmas. In Sect. 3, we establish an solvability and uniqueness of the higher-order fractional neutral stochastic integro-delay differential system. In Sect. 4, we give a sufficient condition of the exponential stability of the system \eqref{eq1}. Finally, in Sect. 5, an example is provided to illustrate the suitability of our results.   
\section {Problem Formulation and Preliminaries}
In this section, we discuss about exponential stability  of higher-order fractional neutral stochastic  differential system driven by Poisson jump is of the form 
 \begin{align}\label{eq1} 
	&^{C}D_{\iota}^{\alpha}\big[x(\iota)+g(\iota, x_{\iota})\big] 	= \ \mathcal{A}x(\iota)+f(\iota, x_{\iota})+G(\iota, x_{\iota})\frac{dw(\iota)}{d\iota}
	 + \int\limits_{\mathcal{Z}}\sigma(\iota, x_{\iota}, u) N(d\iota, du) , \ \iota \in J:=[0, b] \nonumber \\  
	&	x(0) = \ \phi \in L_{p}(\Omega, \mathcal{B}) \ ; \ \frac{d}{d\iota}[x(\iota)+g(\iota, x_{\iota})]|_{\iota=0}  = \eta \in \mathbb{H} ,
\end{align}
\noindent
where $x(\cdot)$ is a stochastic process in a separable Hilbert space $\mathbb{H},$ with the inner product $\langle \cdot , \cdot  \rangle _{\mathbb{H}}$ and the norm $\|\cdot\|_{\mathbb{H}}$. Here, $^{C}D_{\iota}^{\alpha}$ denotes the Caputo fractional derivative of order $1 < \alpha < 2$; $\mathcal{A}: \mathbb{D}(\mathcal{A}) \subset \mathbb{H} \rightarrow \mathbb{H}$ is the infinitesimal generator of an $ C_{0} -$ semigroup of $ \alpha- $ order cosine family $C_{\alpha}(\iota), 0\leq \iota < \infty$ of strongly continuous bounded linear operators associated with sine operator  $S_{\alpha}(\iota), 0\leq \iota < \infty$. The state variable $ x_{\iota}: J \rightarrow \mathbb{H}, x_{\iota}(\theta)= x(\iota+\theta) $ defined on a separable Hilbert space $\mathbb{H}$. Let $ \mathbb{K} $ be the another separable Hilbert and  with the inner product $\langle \cdot , \cdot  \rangle _{\mathbb{K}}$ and the norm $\|\cdot\|_{\mathbb{K}}.$ The initial data $\phi, \eta$ are the $\mathfrak{F}_{0}$-measurable $\mathbb{H}$-valued stochastic process independent of the Brownian motion and $\{w(\iota); \iota \in J\}$ is a standard Wiener process on a real and separable Hilbert space $\mathbb{H}$. In addition $N(d\iota, du)$ denotes the Poisson point process $N(d\iota, du)= \tilde{N}(d\iota, du)-d\iota(\lambda du)$ is the Poisson measure and $N(d\iota, du)$ calls Poisson counting measure correlated with a characteristic measure $\lambda$.  The non-linear maps are $f, g: J \times \mathcal{B}_{r}  \rightarrow \mathbb{H}, G: J \times \mathcal{B}_{r}  \rightarrow L_{Q}^{0}(\mathbb{K},\mathbb{H})$, $\sigma:J \times \mathcal{B}_{r} \rightarrow \mathbb{H}, i=1, 2$ and $\hat{\sigma}:J \times \mathcal{B}_{r}  \times \mathcal{Z} \rightarrow \mathbb{H}$ are appropriate continuous functions to be specified in the sequel. Here, $ L_{Q}^{0} (\mathbb{K},\mathbb{H})$ is the space of all $ Q- $Hilbert Schmidt operators from $ \mathbb{K}$ into $\mathbb{H}.$   	
\subsection{Preliminaries}
In this part, we recollect some basic concepts of fractional calculus (Riemann-Liouville) fractional derivative and, Caputo derivative, stochastic analysis technique, existing lemmas and $ C_{0}- $ semigroups are defined in the sequel. Take $(\Omega, \mathfrak{F},\mathcal{P})$ to be a complete filtered probability space furnished with complete family of right continuous increasing sub $\sigma$-algebras $\{\mathfrak{F}_\iota, \iota\in J\}$ satisfying $\mathfrak{F}_\iota \subset \mathfrak{F}.$ A $\mathbb{H}$-valued random variable is an $\mathfrak{F}_{\iota}$-measurable function $x(\iota):\Omega \rightarrow \mathbb{H}$, and the space  $\mathbb{S}=\{x(\iota,\omega):\Omega \rightarrow \mathbb{H}: \iota\in J\}$ which contains a collection of random variable is called a stochastic process. Let $\gamma_{n}(\iota)(n=1,2, \dots)$ be a sequence of real valued one-dimensional standard Bm  independent of $(\Omega, \mathfrak{F},\mathcal{P})$.
 Set $w(\iota)= \sum_{n = 1}^{\infty} \sqrt{\lambda_{n}} \gamma_{n}(\iota) \zeta_{n}(\iota), \iota \geq 0$, where $\lambda_{n} \geq 0$ are positive real numbers and $\{\zeta_{n}\}(n=1,2, \dots)$ is complete orthonormal basis in $\mathbb{K}$. Let $Q \in L(\mathbb{K}, \mathbb{H})$ be an operator defined by $Q \zeta_{n} = \lambda_{n} \zeta_{n}$ with finite $ Tr(Q) = \sum _{n = 1} ^{\infty} \lambda_{n} < \infty$. Then the above $\mathbb{K}$-valued stochastic process $w(\iota)$ is called as $Q$-Wiener process. Let $\Psi \in L_{Q}^{0}(\mathbb{K}, \mathbb{H}),$ 
 $$ \|\Psi\|^{2}_{Q} = Tr (\Psi Q \Psi^{*}) = \sum \limits_{n= 1}^{\infty} \|\sqrt{\lambda_{n}} \Psi \zeta_{n}\|^{2}.$$ 
 If  $\|\Psi\|_{Q} < \infty$, then $\Psi$ is known as $Q$-Hilbert Schmidt operator.
The articles can be consulted for more information on the concepts and theory of SDEs \cite{r6, r7, r15, r20} and references therein.
Now, we define the space formed by  $\mathfrak{F}_{\iota}-$adapted $ \mathbb{H}-$ valued measurable stochastic process $\{x(\iota), \iota \in J \}$ s.t $x$ is continuous  endowed with the supremum norm
$$ \|x\|^{p} = \Big(\sup\limits_{0 \leq \iota \leq T} \mathbb{E}\|x(\iota)\|^{p}\Big)^{\frac{1}{p}}.$$
Now, define the set
$$\mathcal{B}_{r}= \Big\{x \in \mathcal{C}(J, L_{p}(\Omega, \mathbb{H})), \mathbb{E}\|x(\iota)\|^{p} \leq r,\ \iota \in J \Big\}$$ is a closed ball in $ \mathcal{C}(J, L_{p}(\Omega, \mathbb{H})).$
Then $(\mathcal{C}(\mathbb{H}), \|\cdot\|)$ is Banach space.
\begin{definition} \cite{r6}
	The R-L integral of order $ n - 1 < q < n $, for a continuous function $f: J \rightarrow \mathbb{R}$ 
	$$ I^{q}_{0^{+}} f(\iota) = \frac{1}{\Gamma (q)} \int\limits_{0}^{\iota} (\iota -\upsilon)^{q- 1} f(\upsilon) \iota, \quad \iota > 0, \ q > 0,$$
	provided that the R.H.S is point-wise defined on $J.$
\end{definition}
\begin{definition}\cite{r11}
	The R-L  derivative of order $q$, of a function $ f:J \rightarrow \mathbb{R}$,
	$$ D^{q}_{0^{+}}  f(\iota) =  \frac{1}{\Gamma(n- q)} \Big(\frac{d}{d\iota}\Big)^{n} \int\limits_{0}^{\iota} (\iota- \upsilon)^{ n - q - 1} f(\upsilon) \iota, \quad \iota > 0, \ n-1< q < n,$$
	here $n = [q] + 1$, $[q]$ denotes the integral part of number $q$, provided that the R.H.S is point-wise defined on $J$.
\end{definition}
\begin{definition} \cite{r9} 
	The Caputo derivative  of order $q$ for a function $ f: J \rightarrow \mathbb{R}$,
	$$ ^{C} D^{q}_{0^{+}}  f(\iota) =  \frac{1}{\Gamma(n- q)} \int\limits_{0}^{\iota} (\iota- \upsilon)^{n- q - 1}  f^{(n)}(\upsilon) \iota, \ n- 1 < q < n ,$$
	where $(n)$ denotes the $n^{th}$ derivative, provided that the R.H.S is point-wise defined on $J$. 	
\end{definition}
\begin{definition}\cite{r13, r17}
	A one parameter family $(C_{\alpha}(\iota))_{\iota\geq 0}, \alpha \in (1, 2]$ is called the solution operator (or) strongly continuous $\alpha$-order fractional cosine family for the given Cauchy problem \eqref{eq1} and $\mathcal{A}$ is the infinitesimal generator of $\{C_{\alpha}(\iota)\}_{\iota\geq 0}$  if the following conditions are satisfied  
	\begin{enumerate}
		\item $C_{\alpha}(\iota)$ is strongly continuous for $\iota\geq 0$ and $C_{\alpha}(0)= I$, where $I$ is an identity operator.
		\item $C_{\alpha}(\iota)\mathbb{D}(\mathcal{A})\subset \mathbb{D}(\mathcal{A}) $ and $\mathcal{A}C_{\alpha}(\iota)\eta= C_{\alpha}(\iota)\mathcal{A} \eta \ \eta\in \mathbb{D}(\mathcal{A}, \iota \geq 0$;
		\item $C_{\alpha}(\iota) \eta$ is solution of $x(\iota) = \eta+ \frac{1}{\Gamma s} \int\limits_{0}^{\upsilon}(\iota-\upsilon)^{\alpha-1}\mathcal{A}x(\upsilon)d\upsilon,  \forall \eta\in \mathbb{D}(\mathcal{A}).$
	\end{enumerate} 
\end{definition}
\begin{definition}\cite{r18, r19}
The fractional sine family $S_{\alpha}(\iota): [0, \infty) \rightarrow \mathbb{R}^{+}$ associated with $C_{\alpha}(\iota)$ is defined by
\begin{align*}
S_{\alpha}(\iota) &= \int\limits_{0}^{\iota} C_{\alpha}(\upsilon)d\upsilon, \iota \leq 0.
\end{align*}
\end{definition}
\begin{definition}\cite{r12, r14}
The fractional R-L family $P_{\alpha}:[0, \infty) \rightarrow \mathbb{R}^{+}$ associated with $C_{\alpha}(\iota)$ is defined by
$$P_{\alpha}(\iota)=J^{\alpha-1} C_{\alpha}(\iota).$$
\end{definition}
\begin{definition}\label{def2.7} \cite{r18}
For any $x \in \mathbb{H}, \ 0 < \gamma < 1$ and $\mu\in (0,1]$, we have $\mathcal{A}S_{\alpha}(\iota)x = \mathcal{A}^{1-\gamma}S_{\alpha}(\iota)\mathcal{A}^{\gamma} x$ and
$$\|\mathcal{A}^{1-\gamma}S_{\alpha}(\iota)\|\leq \frac{\alpha c_{\mu}}{\iota^{\alpha \mu}}, \ \forall \quad \iota \in J.$$
\end{definition}
\begin{lemma}\cite{r13} \label{2.6}
For any fixed $\iota\geq 0$ and for any $x \in \mathbb{H},$ the following estimates are true
$$\|C_{\alpha}(\iota)\|\leq M_{c}; \ \|S_{\alpha}(\iota)\|\leq M\iota.$$
\end{lemma}
\begin{lemma}\label{lm2} \cite{xm, r5} 
	For any $p \geq 1$ and for an arbitrary $L_{Q}^{0} (\mathbb{K}, \mathbb{H})$-valued predictable process $G(\cdot)$ such that  
	\begin{align*}    
	\sup \limits_{\iota \in J} \mathbb{E} \Big \| \int \limits_{0}^{\iota} G(\upsilon) dw(\upsilon) \Big \|^{p} &\leq  \tilde{C}_{p}\left(\int\limits _{0} ^{\iota} (\mathbb{E}\|G(\upsilon)\|^{p}_{L_{Q}^{0}} )^{\frac{2}{p}}d\upsilon  \right)^{\frac{p}{2}},   \ \iota \in J, 
	\end{align*}
where $\tilde{C}_{p} = \Big(\frac{p(p -1)}{2}\Big)^{\frac{p}{2}}.$
\end{lemma}
\begin{lemma}\label{2.8}\cite{k}
For any $p \geq 2$, there exists $\tilde{\mathcal{C}_{p}} > 0$ such that
\begin{align*}
\sup \limits_{0 \leq \iota \leq \iota}  \mathbb{E}\Big [ \|\int \limits_{0}^{\iota} \int \limits_{\mathbb{Z}} \sigma (\iota, x) \tilde{N}(d \iota, dx) \|\Big] ^{p}  \leq \tilde{\mathcal{C}_{p}} \left\{\mathbb{E}\Big[\Big( \int \limits _{0}^{\iota} \int\limits _{\mathbb{Z}} \| \sigma(\iota,x) \|^{2} \lambda dx d\iota\Big)^{\frac{p}{2}}\Big]\right.\\
\left.+ \mathbb{E}\Big[ \int \limits _{0}^{\iota} \int\limits _{\mathbb{Z}} \| \sigma(\iota,x) \|^{2p} \lambda dx d\iota\Big]^{\frac{1}{2}}\right\}. 
\end{align*}
\end{lemma}

\begin{definition}
An $X$-valued stochastic process $x(\iota)\ (\iota \in J),$ is said to be a mild solution for system \eqref{eq1}, provided 
\begin{enumerate}
\item $x(\iota)$ is fitting a $\mathfrak{F}_{\iota}\ (\iota \geq 0)$ has a c$\grave{a}$dl$\grave{a}$g path  $\iota \geq 0$ a.s.
\item For $\iota \in J$, the following integral equation is satisfied:  
\end{enumerate}
	\begin{equation*}
x(\iota) = \left\{\begin{array}{ll}
\phi(\iota),  \quad\quad \iota \in (-\infty, 0];\nonumber\\
C_{\alpha}(\iota)[\phi+g(0, \phi)] +S_{\alpha}(\iota)\eta -g(\iota, x_{\iota})- \int\limits_{0}^{\iota}\mathcal{A}S_{\alpha}(\iota-\upsilon)g(\upsilon, x_{\upsilon})d\upsilon\nonumber\\
+ \int\limits_{0}^{\iota}S_{\alpha}(\iota-\upsilon) f(\upsilon, x_{\upsilon})d\upsilon 
+ \int\limits_{0}^{\iota}S_{\alpha}(\iota-\upsilon) G(\upsilon, x_{\upsilon})dw(\upsilon)\nonumber \\
+ \int\limits_{0}^{\iota}\int\limits_{\mathcal{Z}}S_{\alpha}(\iota-\upsilon) \sigma(\upsilon, x_{\upsilon}, u)N(d\upsilon, du),
   \quad \iota\in J.
		\end{array}\right. 
		\end{equation*}
	\end{definition}
\begin{definition}\label{def2.9}
Suppose that there exists mappings are $ \Gamma_{1}: J \times \mathcal{B} \rightarrow \mathbb{H};$ $ \Gamma_{2}: J \times \mathcal{B} \rightarrow  \mathbb{H}; \Gamma_{3}: J \times \mathcal{B}\rightarrow L_{Q}^{0}(\mathbb{K}, \mathbb{H})$ and $\Gamma_{4}: J \times \mathcal{B} \times \mathcal{Z}\rightarrow \mathbb{H}$ be real valued functions on a bounded linear operator. Then there exist positive constants $ \hat{a}_{i}, i= 1, 2, 3, 4, 5 $ s.t for any $ x, y \in \mathbb{H},$ we have
 \begin{align}\label{eq3}
&(i) \ \mathbb{E}\Big\|f\Big(\iota, x(\iota)+ y(\iota)+\int\limits_{0}^{\iota}S_{\alpha}(\iota-\upsilon)\Gamma_{1} (\upsilon, x(\upsilon)) y(\upsilon)d\upsilon+\int\limits_{0}^{\iota}S_{\alpha}(\iota-\upsilon)\Gamma_{2} (\upsilon, x(\upsilon)) y(\upsilon)d\upsilon\nonumber\\
&+\int\limits_{0}^{\iota}S_{\alpha}(\iota-\upsilon)\Gamma_{3} (\upsilon, x(\upsilon)) y(\upsilon)dw(\upsilon)+\int\limits_{0}^{\iota}\int\limits_{\mathcal{Z}}S_{\alpha}(\iota-\upsilon)\Gamma_{4} (\upsilon, x(\upsilon), u) y(\upsilon)\lambda du\Big)-f(\iota, x(\iota))\nonumber\\
&- \Gamma_{1} (\iota, x(\iota))y(\iota)\Big\|^{p} \leq \hat{a}_{2} \mathbb{E}\|y(\iota)\|^{p}.\nonumber
\end{align}
\begin{align}
&(ii) \ \mathbb{E}\Big\|\mathcal{A}^{\gamma}g\Big(\iota, x(\iota)+ y(\iota)+\int\limits_{0}^{\iota}S_{\alpha}(\iota-\upsilon)\Gamma_{1} (\upsilon, x(\upsilon)) y(\upsilon)d\upsilon+\int\limits_{0}^{\iota}S_{\alpha}(\iota-\upsilon)\Gamma_{2} (\upsilon, x(\upsilon)) y(\upsilon)d\upsilon\nonumber\\
&+\int\limits_{0}^{\iota}S_{\alpha}(\iota-\upsilon)\Gamma_{3} (\upsilon, x(\upsilon)) y(\upsilon)dw(\upsilon)+\int\limits_{0}^{\iota}\int\limits_{\mathcal{Z}}S_{\alpha}(\iota-\upsilon)\Gamma_{4} (\upsilon, x(\upsilon), u) y(\upsilon)\lambda du\Big)\nonumber\\
&-\mathcal{A}^{\gamma}g(\iota, x(\iota))- \Gamma_{1} (\iota, x(\iota))y(\iota)\Big\|^{p} \leq \hat{a}_{1} \mathbb{E}\|y(\iota)\|^{p}.\nonumber\\
&(iii)\ \mathbb{E}\Big\|G\Big(\iota, x(\iota)+ y(\iota)+\int\limits_{0}^{\iota}S_{\alpha}(\iota-\upsilon)\Gamma_{1} (\upsilon, x(\upsilon)) y(\upsilon)d\upsilon+\int\limits_{0}^{\iota}S_{\alpha}(\iota-\upsilon)\Gamma_{2} (\upsilon, x(\upsilon)) y(\upsilon)d\upsilon\nonumber\\
&+\int\limits_{0}^{\iota}S_{\alpha}(\iota-\upsilon)\Gamma_{3} (\upsilon, x(\upsilon)) y(\upsilon)dw(\upsilon)+\int\limits_{0}^{\iota}\int\limits_{\mathcal{Z}}S_{\alpha}(\iota-\upsilon)\Gamma_{4} (\upsilon, x(\upsilon), u) y(\upsilon)\lambda du\Big)-G(\iota, x(\iota))\nonumber\\
&- \Gamma_{3} (\iota, x(\iota))y(\iota)\Big\|^{p} \leq \hat{a}_{3} \mathbb{E}\|y(\iota)\|^{p}.\nonumber\\
&(iv)\ \Big[\int\limits_{\mathcal{Z}}\mathbb{E}\Big\|\sigma\Big(\iota, x(\iota)+ y(\iota)+\int\limits_{0}^{\iota}S_{\alpha}(\iota-\upsilon)\Gamma_{1} (\upsilon, x(\upsilon)) y(\upsilon)d\upsilon+\int\limits_{0}^{\iota}S_{\alpha}(\iota-\upsilon)\Gamma_{2} (\upsilon, x(\upsilon)) y(\upsilon)d\upsilon\nonumber\\
&+\int\limits_{0}^{\iota}S_{\alpha}(\iota-\upsilon)\Gamma_{3} (\upsilon, x(\upsilon)) y(\upsilon)dw(\upsilon)+\int\limits_{0}^{\iota}\int\limits_{\mathcal{Z}}S_{\alpha}(\iota-\upsilon)\Gamma_{4} (\upsilon, x(\upsilon), u) y(\upsilon)\lambda du\Big)\lambda du\nonumber\\
&-\int\limits_{\mathcal{Z}}\sigma(\iota, x(\iota), u)N(d\iota, du)- \Gamma_{4} (\iota, x(\iota))y(\iota)\Big\|^{2}\lambda du\Big]^{\frac{p}{2}} \leq \hat{a}_{4} \mathbb{E}\|y(\iota)\|^{p}.\nonumber\\
&(v)\ \Big[\int\limits_{\mathcal{Z}}\mathbb{E}\Big\|\sigma\Big(\iota, x(\iota)+ y(\iota)+\int\limits_{0}^{\iota}S_{\alpha}(\iota-\upsilon)\Gamma_{1} (\upsilon, x(\upsilon)) y(\upsilon)d\upsilon+\int\limits_{0}^{\iota}S_{\alpha}(\iota-\upsilon)\Gamma_{2} (\upsilon, x(\upsilon)) y(\upsilon)d\upsilon\nonumber\\
&+\int\limits_{0}^{\iota}S_{\alpha}(\iota-\upsilon)\Gamma_{3} (\upsilon, x(\upsilon)) y(\upsilon)dw(\upsilon)+\int\limits_{0}^{\iota}\int\limits_{\mathcal{Z}}S_{\alpha}(\iota-\upsilon)\Gamma_{4} (\upsilon, x(\upsilon), u) y(\upsilon)\lambda du\Big)\lambda du\nonumber\\
&-\int\limits_{\mathcal{Z}}\sigma(\iota, x(\iota), u)N(d\iota, du)- \Gamma_{4} (\iota, x(\iota))y(\iota)\Big\|^{2p}\lambda du\Big]^{\frac{1}{2}} \leq \hat{a}_{5} \mathbb{E}\|y(\iota)\|^{p}.
\end{align}
Then, $\Gamma_{1}, \Gamma_{2}, \Gamma_{3}, \Gamma_{4}$ are integral contractors for the non-linear functions $f, g, G, \sigma$, respectively.
\end{definition}		
\begin{remark}
If $ \Gamma_{1}\equiv \Gamma_{2}\equiv, \Gamma_{3} \equiv \Gamma_{4} \equiv 0,$ then \eqref{eq3} reduces to Lipschitz conditions (i.e)
\begin{align*}
\mathbb{E}\Big\|f(\cdot, x(\cdot)+ y(\cdot)-f(\cdot, x(\cdot)))\Big\|^{p} & \leq \hat{a}_{2} \mathbb{E}\|y(\cdot)\|^{p},\\
\mathbb{E}\Big\|g(\cdot, x(\cdot) + y(\cdot)-g(\cdot, x(\cdot)))\Big\|^{p} &\leq \hat{a}_{1} \mathbb{E}\|y(\cdot)\|^{p},\\
\mathbb{E}\Big\|G(\cdot, x(\cdot) + y(\cdot)-G(\cdot, x(\cdot)))\Big\|^{p} &\leq \hat{a}_{3} \mathbb{E}\|y(\cdot)\|^{p},
\end{align*}
\begin{align*}
\mathbb{E}\Big\|\int\limits_{\mathcal{Z}}\Big[\sigma(\cdot, x(\cdot) + y(\cdot), u)-\sigma(\cdot, x(\cdot), u)\Big]\lambda du\Big\|^{p} &\leq \hat{a}_{4} \mathbb{E}\|y(\cdot)\|^{p}.
\end{align*}
\end{remark}
\begin{remark}
It should be remarked here that the Lipschitz conditions leads to a unique solution of \eqref{eq1}, but the condition given in \eqref{eq3} may not give the uniqueness of the solution of \eqref{eq1}. The uniqueness of the solution of the given system is ensured by the regularity of the integral contractor. 
 \end{remark}		
\begin{definition}\label{def2.11}
A bounded stochastic integral contractor is said to be regular, if the integral equation
\begin{align} 
&y(\iota)+\int\limits_{0}^{\iota}S_{\alpha}(\iota-\upsilon)\Gamma_{1} (\upsilon, x(\upsilon))y(\upsilon) +\int\limits_{0}^{\iota}S_{\alpha}(\iota-\upsilon)\Gamma_{2} (\upsilon, x(\upsilon))y(\upsilon)\nonumber\\
&+\int\limits_{0}^{\iota}S_{\alpha}(\iota-\upsilon)\Gamma_{3} (\upsilon, x(\upsilon))y(\upsilon)dw(\upsilon)+\int\limits_{0}^{\iota}\int\limits_{\mathcal{Z}}S_{\alpha}(\iota-\upsilon)\Gamma_{4} (\upsilon, x(\upsilon), u)y(\upsilon)\lambda du= A(\iota)
\end{align}
has a solution  $ y \in \mathbb{H} $ for any $ x, A \in \mathbb{H}.$\\ Let us assume that $ \mathbb{E}\|\Gamma_{1}(\iota, x(\iota))\|^{p}\leq c_{1},\ \mathbb{E}\|\Gamma_{2}(\iota, x(\iota))\|^{p}\leq c_{2} $, $\mathbb{E}\|\Gamma_{3}(\iota, x(\iota))\|^{p}\leq c_{3} $, and $\mathbb{E}\|\Gamma_{4}(\iota, x(\iota))\|^{\frac{p}{2}}\leq c_{4} $, $\mathbb{E}\|\Gamma_{4}(\iota, x(\iota))\|^{p}\leq \hat{c}_{4} $ for all $ \iota \in J,\ x \in \mathbb{H}.$
\end{definition}		
\section{Existence and Uniqueness of the Proposed System}
In this section, we show that the existence and uniqueness of mild solution for the proposed system \eqref{eq1}. After that proving the existence results for the mild solution, the following hypotheses are necessary to prove the main results:
\begin{itemize}
	\item [($H_1$)] The non-linear continuous functions $f, g, G$ and 
$\sigma$ have  regular integral contractors $ \Gamma_{1}, \Gamma_{2}, \Gamma_{3}$, and $\Gamma_{4}, respectively.$  
\end{itemize}
\begin{theorem}\label{PIDthm1}
If the	assumption ($H_{1}$) holds, then the fractional neutral higher-order stochastic differential system has a unique mild solution defined on $ J, $ provided
\begin{align*}
&5^{p-1} \Big\{\|\mathcal{A}^{-\gamma}\|^{p} \ \hat{a}_{1}^{p}+\alpha^{p}\ c_{\mu}^{p}\ \hat{a}_{1}^{p}\ \Big[\frac{\iota^{p\alpha \mu}}{p\alpha \mu}\Big]^{p}+k(p)\ M^{p}\Big\{ \Big(\frac{\iota^{3}}{3}\Big)^{\frac{p}{2}}\hat{a}_{4}^{p}+ \Big(\frac{\iota^{2p+1}}{2p+1}\Big)^{\frac{1}{2}}\hat{a}_{5}^{p}\Big\} \nonumber\\
	&\quad+M^{p} \ \iota^{p-1} \  \hat{a}_{2}^{p}+ C_{p}\  M^{P}\ \iota^{\frac{p}{2}-1} \ \hat{a}_{3}^{p}\Big\}< 1.
\end{align*} 
\end{theorem}
\noindent\textbf{Proof}: 
\textbf{Existence of Mild Solution}\\
We consider the operator $ \Phi: \mathcal{B}_{r} \rightarrow \mathcal{B}_{r}$ defined by
	\begin{equation*}
(\Phi x)(\iota) = \left\{\begin{array}{ll}
\phi(\iota),  \quad\quad \iota \in (-\infty, 0]\nonumber\\
C_{\alpha}(\iota)[\phi+g(0, \phi)] +S_{\alpha}(\iota)\eta -g(\iota, x_{\iota})- \int\limits_{0}^{\iota}\mathcal{A}S_{\alpha}(\iota-\upsilon)g(\upsilon, x_{\upsilon})d\upsilon\nonumber\\
+ \int\limits_{0}^{\iota}S_{\alpha}(\iota-\upsilon) f(\upsilon, x_{\upsilon})d\upsilon 
+ \int\limits_{0}^{\iota}S_{\alpha}(\iota-\upsilon) G(\upsilon, x_{\upsilon})dw(\upsilon)\nonumber \\
+ \int\limits_{0}^{\iota}\int\limits_{\mathcal{Z}}S_{\alpha}(\iota-\upsilon) \sigma(\upsilon, x_{\upsilon}, u)N(d\upsilon, du),
   \quad \iota\in J.
		\end{array}\right. 
		\end{equation*}
In order to prove that the existence and unique result for the mild solution of \eqref{eq1}, it is enough to show that the operator $ \Phi $ has a unique mild solution. By using successive approximation technique, we study the following cases: 

\noindent
Case (i): Consider the following two sequences $ \{x_{n}\}_{n=1}^{\infty} $ and $ \{y_{n}\}_{n=1}^{\infty} $
in $ \mathbb{H}.$
	\begin{align}
	x_{0}(\iota) &= \ C_{\alpha}(\iota)[\phi+g(0, \phi)]+S_{\alpha}(\iota)\eta,\nonumber\\
	y_{n}(\iota)&=\ x_{n}(\iota)-g(\iota, x_{\iota}^{n})- \int\limits_{0}^{\iota}\mathcal{A}S_{\alpha}(\iota-\upsilon)g(\upsilon, x_{\upsilon}^{n})d\upsilon+ \int\limits_{0}^{\iota}S_{\alpha}(\iota-\upsilon) f(\upsilon, x_{\upsilon}^{n})d\upsilon \nonumber\\
	&\quad
	+ \int\limits_{0}^{\iota}S_{\alpha}(\iota-\upsilon) G(\upsilon, x_{\upsilon}^{n})dw(\upsilon)+ \int\limits_{0}^{\iota}\int\limits_{\mathcal{Z}}S_{\alpha}(\iota-\upsilon) \sigma(\upsilon, x_{\upsilon}^{n}, u)N(d\upsilon, du).
	\end{align}
	\begin{align*}
	x_{n+1}(\iota)&=\ x_{n}(\iota)- \Big[y_{n}(\iota)+ \int\limits_{0}^{\iota}S_{\alpha}(\iota-\upsilon)\Gamma_{1}(\upsilon, x_{n}(\upsilon)) y_{n}(\upsilon)+ \int\limits_{0}^{\iota}S_{\alpha}(\iota-\upsilon)\Gamma_{2}(\upsilon, x_{n}(\upsilon)) y_{n}(\upsilon)\nonumber\\
	&\quad+ \int\limits_{0}^{\iota}S_{\alpha}(\iota-\upsilon)\Gamma_{3}(\upsilon, x_{n}(\upsilon)) y_{n}(\upsilon)dw(\upsilon)+\int\limits_{0}^{\iota}\int\limits_{\mathcal{Z}}S_{\alpha}(\iota-\upsilon)\Gamma_{4} (\upsilon, x_{n}(\upsilon), u) y_{n}(\upsilon)\lambda du\Big]\nonumber\\
	&=\ -g(\iota,x^{n}_{\iota})-\int\limits_{0}^{\iota}S_{\alpha}(\iota-\upsilon)\Gamma_{1}(\upsilon, x_{n}(\upsilon)) y_{n}(\upsilon)d\upsilon- \int\limits_{0}^{\iota}\mathcal{A}S_{\alpha}(\iota-\upsilon)g(\upsilon, x_{\upsilon}^{n})d\upsilon\nonumber\\
		&\quad-\int\limits_{0}^{\iota}S_{\alpha}(\iota-\upsilon)\Gamma_{1}(\upsilon, x_{n}(\upsilon)) y_{n}(\upsilon)d\upsilon+ \int\limits_{0}^{\iota}S_{\alpha}(\iota-\upsilon) f(\upsilon,x_{\upsilon}^{n})d\upsilon-\int\limits_{0}^{\iota}S_{\alpha}(\iota-\upsilon)\Gamma_{2}(\upsilon, x_{n}(\upsilon)) y_{n}(\upsilon)d\upsilon \nonumber\\
	&\quad+ \int\limits_{0}^{\iota}S_{\alpha}(\iota-\upsilon) G(\upsilon, x_{\upsilon}^{n})dw(\upsilon)-\int\limits_{0}^{\iota}S_{\alpha}(\iota-\upsilon)\Gamma_{2}(\upsilon, x_{n}(\upsilon)) y_{n}(\upsilon)dw(\upsilon)\nonumber\\
	&\quad+ \int\limits_{0}^{\iota}\int\limits_{\mathcal{Z}}S_{\alpha}(\iota-\upsilon) \sigma(\upsilon, x_{\upsilon}^{n}, u)N(d\upsilon, du)-\int\limits_{0}^{\iota}\int\limits_{\mathcal{Z}}S_{\alpha}(\iota-\upsilon)\Gamma_{4} (\upsilon, x_{n}(\upsilon), u) y_{n}(\upsilon)\lambda du+x_{0}(\iota).	
	\end{align*}
Now, we have  
\begin{align*}
y_{n+1}(\iota)&=\ x_{n+1}(\iota)-g(\iota, x_{\iota}^{n+1})- \int\limits_{0}^{\iota}\mathcal{A}S_{\alpha}(\iota-\upsilon)g(\upsilon, x_{\upsilon}^{n+1})d\upsilon+ \int\limits_{0}^{\iota}S_{\alpha}(\iota-\upsilon) f(\upsilon, x_{\upsilon}^{n+1})d\upsilon \nonumber\\
	&\quad
	+ \int\limits_{0}^{\iota}S_{\alpha}(\iota-\upsilon) G(\upsilon, x_{\upsilon}^{n+1})dw(\upsilon)+ \int\limits_{0}^{\iota}\int\limits_{\mathcal{Z}}S_{\alpha}(\iota-\upsilon) \sigma(\upsilon, x_{\upsilon}^{n+1}, u)N(d\upsilon, du)-x_{0}(\iota).
	\end{align*}
	Now, we have to substitute the value of $ x_{n+1} $ in the above inequality, we get the following
	\begin{align*}
	&\quad= -g(\iota,x^{n}_{\iota})-\int\limits_{0}^{\iota}S_{\alpha}(\iota-\upsilon)\Gamma_{1}(\upsilon, x_{\upsilon}^{n}) y_{n}(\upsilon)d\upsilon- \int\limits_{0}^{\iota}\mathcal{A}S_{\alpha}(\iota-\upsilon)g(\upsilon, x_{\upsilon}^{n})d\upsilon\nonumber\\
&\quad-\int\limits_{0}^{\iota}S_{\alpha}(\iota-\upsilon)\Gamma_{1}(\upsilon, x_{\upsilon}^{n}) y_{n}(\upsilon)d\upsilon+ \int\limits_{0}^{\iota}S_{\alpha}(\iota-\upsilon) f(\upsilon,x_{\upsilon}^{n})d\upsilon-\int\limits_{0}^{\iota}S_{\alpha}(\iota-\upsilon)\Gamma_{2}(\upsilon, x_{\upsilon}^{n}) y_{n}(\upsilon)d\upsilon \nonumber\\
&\quad+ \int\limits_{0}^{\iota}S_{\alpha}(\iota-\upsilon) G(\upsilon, x_{\upsilon}^{n})dw(\upsilon)-\int\limits_{0}^{\iota}S_{\alpha}(\iota-\upsilon)\Gamma_{2}(\upsilon, x_{\upsilon}^{n}) y_{n}(\upsilon)dw(\upsilon)\nonumber\\
&\quad+ \int\limits_{0}^{\iota}\int\limits_{\mathcal{Z}}S_{\alpha}(\iota-\upsilon) \sigma(\upsilon, x_{\upsilon}^{n}, u)N(d\upsilon, du)-\int\limits_{0}^{\iota}\int\limits_{\mathcal{Z}}S_{\alpha}(\iota-\upsilon)\Gamma_{4} (\upsilon, x_{\upsilon}^{n}, u) y_{n}(\upsilon)\lambda du\\	
& -g\Big(\iota, x^{n}_{\iota}- y^{n}_{\iota}-\int\limits_{0}^{\iota}S_{\alpha}(\iota-\upsilon)\Gamma_{1}(\upsilon, x_{\upsilon}^{n}) y_{n}(\upsilon)d\upsilon-\int\limits_{0}^{\iota}S_{\alpha}(\iota-\upsilon)\Gamma_{2}(\upsilon, x_{\upsilon}^{n}) y_{n}(\upsilon)d\upsilon\nonumber\\
	&\quad- \int\limits_{0}^{\iota}S_{\alpha}(\iota-\upsilon)\Gamma_{3}(\upsilon, x_{\upsilon}^{n}) y_{n}(\upsilon)dw(\upsilon)-\int\limits_{0}^{\iota}\int\limits_{\mathcal{Z}}S_{\alpha}(\iota-\upsilon)\Gamma_{4} (\upsilon, x_{\upsilon}^{n}, u) y_{n}(\upsilon)\lambda du\Big)\nonumber\\
&\quad-	\int\limits_{0}^{\iota}\mathcal{A}S_{\alpha}(\iota-\upsilon)g\Big(\upsilon, x^{n}_{\upsilon}-y^{n}_{\upsilon}- \int\limits_{0}^{\kappa}S_{\alpha}(\kappa-\upsilon)\Gamma_{1}(\upsilon, x_{\upsilon}^{n}) y_{n}(\upsilon)d\upsilon- \int\limits_{0}^{\kappa}S_{\alpha}(\kappa-\upsilon)\Gamma_{2}(\upsilon, x_{\upsilon}^{n}) y_{n}(\upsilon)d\upsilon\nonumber\\
	&\quad- \int\limits_{0}^{\kappa}S_{\alpha}(\kappa-\upsilon)\Gamma_{3}(\upsilon, x_{\upsilon}^{n}) y_{n}(\upsilon)dw(\upsilon)-\int\limits_{0}^{\kappa}\int\limits_{\mathcal{Z}}S_{\alpha}(\kappa-\upsilon)\Gamma_{4} (\upsilon, x_{\upsilon}^{n}, u) y_{n}(\upsilon)N(d\upsilon, du)\Big)d\upsilon\nonumber\\
	&\quad+\int\limits_{0}^{\iota}S_{\alpha}(\iota-\upsilon)f\Big(\upsilon, x^{n}_{\upsilon}-y^{n}_{\upsilon}- \int\limits_{0}^{\kappa}S_{\alpha}(\kappa-\upsilon)\Gamma_{1}(\upsilon, x_{\upsilon}^{n}) y_{n}(\upsilon)d\upsilon- \int\limits_{0}^{\kappa}S_{\alpha}(\kappa-\upsilon)\Gamma_{2}(\upsilon, x_{\upsilon}^{n}) y_{n}(\upsilon)d\upsilon\nonumber\\
		&\quad- \int\limits_{0}^{\kappa}S_{\alpha}(\kappa-\upsilon)\Gamma_{3}(\upsilon, x_{\upsilon}^{n}) y_{n}(\upsilon)dw(\upsilon)-\int\limits_{0}^{\kappa}\int\limits_{\mathcal{Z}}S_{\alpha}(\kappa-\upsilon)\Gamma_{4} (\upsilon, x_{\upsilon}^{n}, u) y_{n}(\upsilon)N(d\upsilon, du)\Big)d\upsilon\nonumber\\
		&\quad+\int\limits_{0}^{\iota}S_{\alpha}(\iota-\upsilon)G\Big(\upsilon, x^{n}_{\upsilon}-y^{n}_{\upsilon}- \int\limits_{0}^{\kappa}S_{\alpha}(\kappa-\upsilon)\Gamma_{1}(\upsilon, x_{\upsilon}^{n}) y_{n}(\upsilon)d\upsilon- \int\limits_{0}^{\kappa}S_{\alpha}(\kappa-\upsilon)\Gamma_{2}(\upsilon, x_{\upsilon}^{n}) y_{n}(\upsilon)d\upsilon\nonumber
			\end{align*}
		\begin{align}\label{eq5}
				&\quad- \int\limits_{0}^{\kappa}S_{\alpha}(\kappa-\upsilon)\Gamma_{3}(\upsilon, x_{\upsilon}^{n}) y_{n}(\upsilon)dw(\upsilon)-\int\limits_{0}^{\kappa}\int\limits_{\mathcal{Z}}S_{\alpha}(\kappa-\upsilon)\Gamma_{4} (\upsilon, x_{\upsilon}^{n}, u) y_{n}(\upsilon)N(d\upsilon, du)\Big)dw(\upsilon)\nonumber\\
				&\quad+\int\limits_{0}^{\iota}\int\limits_{\mathcal{Z}}S_{\alpha}(\iota-\upsilon)\sigma\Big(\upsilon, x^{n}_{\upsilon}-y^{n}_{\upsilon}- \int\limits_{0}^{\kappa}S_{\alpha}(\kappa-\upsilon)\Gamma_{1}(\upsilon, x_{\upsilon}^{n}) y_{n}(\upsilon)d\upsilon- \int\limits_{0}^{\kappa}S_{\alpha}(\kappa-\upsilon)\Gamma_{2}(\upsilon, x_{\upsilon}^{n}) y_{n}(\upsilon)d\upsilon\nonumber\\
						&\quad- \int\limits_{0}^{\kappa}S_{\alpha}(\kappa-\upsilon)\Gamma_{3}(\upsilon, x_{\upsilon}^{n}) y_{n}(\upsilon)dw(\upsilon)-\int\limits_{0}^{\kappa}\int\limits_{\mathcal{Z}}S_{\alpha}(\kappa-\upsilon)\Gamma_{4} (\upsilon, x_{\upsilon}^{n}, u) y_{n}(\upsilon)N(d\upsilon, du)\Big)N(d\upsilon, du). 
				\end{align}
Substituting $ x_{n}=x $ and $ y_{n}=-y $ in equation \eqref{eq5}, and by using Definition \ref{def2.11}, we get 
\begin{align}\label{eq6}
y_{n+1}(\iota)& =\ -g(\iota,x_{\iota})+\int\limits_{0}^{\iota}S_{\alpha}(\iota-\upsilon)\Gamma_{1}(\upsilon, x_{\upsilon}) y(\upsilon)d\upsilon- \int\limits_{0}^{\iota}\mathcal{A}S_{\alpha}(\iota-\upsilon)g(\upsilon, x_{\upsilon})d\upsilon\nonumber\\
&\quad+\int\limits_{0}^{\iota}S_{\alpha}(\iota-\upsilon)\Gamma_{1}(\upsilon, x_{\upsilon}) y(\upsilon)d\upsilon+ \int\limits_{0}^{\iota}S_{\alpha}(\iota-\upsilon) f(\upsilon,x_{\upsilon})d\upsilon+\int\limits_{0}^{\iota}S_{\alpha}(\iota-\upsilon)\Gamma_{2}(\upsilon, x_{\upsilon}) y(\upsilon)d\upsilon \nonumber\\
&\quad+ \int\limits_{0}^{\iota}S_{\alpha}(\iota-\upsilon) G(\upsilon, x_{\upsilon})dw(\upsilon)+\int\limits_{0}^{\iota}S_{\alpha}(\iota-\upsilon)\Gamma_{2}(\upsilon, x_{\upsilon}) y(\upsilon)dw(\upsilon)\nonumber\\
&\quad+ \int\limits_{0}^{\iota}\int\limits_{\mathcal{Z}}S_{\alpha}(\iota-\upsilon) \sigma(\upsilon, x_{\upsilon}, u)N(d\upsilon, du)+\int\limits_{0}^{\iota}\int\limits_{\mathcal{Z}}S_{\alpha}(\iota-\upsilon)\Gamma_{4} (\upsilon, x_{\upsilon}, u) y(\upsilon)\lambda du\nonumber\\
	& -g\Big(\iota, x_{\iota}+ y_{\iota}+\int\limits_{0}^{\iota}S_{\alpha}(\iota-\upsilon)\Gamma_{1}(\upsilon, x_{\upsilon}) y(\upsilon)d\upsilon+\int\limits_{0}^{\iota}S_{\alpha}(\iota-\upsilon)\Gamma_{2}(\upsilon, x_{\upsilon}) y(\upsilon)d\upsilon\nonumber\\
	&\quad+ \int\limits_{0}^{\iota}S_{\alpha}(\iota-\upsilon)\Gamma_{3}(\upsilon, x_{\upsilon}) y(\upsilon)dw(\upsilon)+\int\limits_{0}^{\iota}\int\limits_{\mathcal{Z}}S_{\alpha}(\iota-\upsilon)\Gamma_{4} (\upsilon, x_{\upsilon}, u) y(\upsilon)\lambda du\Big)\nonumber\\
&\quad+	\int\limits_{0}^{\iota}\mathcal{A}S_{\alpha}(\iota-\upsilon)g\Big(\upsilon, x_{\upsilon}+y_{\upsilon}+ \int\limits_{0}^{\kappa}S_{\alpha}(\kappa-\upsilon)\Gamma_{1}(\upsilon, x_{\upsilon}) y(\upsilon)d\upsilon+ \int\limits_{0}^{\kappa}S_{\alpha}(\kappa-\upsilon)\Gamma_{2}(\upsilon, x_{\upsilon}) y(\upsilon)d\upsilon\nonumber\\
	&\quad+ \int\limits_{0}^{\kappa}S_{\alpha}(\kappa-\upsilon)\Gamma_{3}(\upsilon, x_{\upsilon}) y(\upsilon)dw(\upsilon)+\int\limits_{0}^{\kappa}\int\limits_{\mathcal{Z}}S_{\alpha}(\kappa-\upsilon)\Gamma_{4} (\upsilon, x_{\upsilon}, u) y(\upsilon)N(d\upsilon, du)\Big)d\upsilon\nonumber\\
	&\quad+\int\limits_{0}^{\iota}S_{\alpha}(\iota-\upsilon)f\Big(\upsilon, x_{\upsilon}+y_{\upsilon}+ \int\limits_{0}^{\kappa}S_{\alpha}(\kappa-\upsilon)\Gamma_{1}(\upsilon, x_{\upsilon}) y(\upsilon)d\upsilon+ \int\limits_{0}^{\kappa}S_{\alpha}(\kappa-\upsilon)\Gamma_{2}(\upsilon, x_{\upsilon}) y(\upsilon)d\upsilon\nonumber\\
		&\quad+ \int\limits_{0}^{\kappa}S_{\alpha}(\kappa-\upsilon)\Gamma_{3}(\upsilon, x_{\upsilon}) y(\upsilon)dw(\upsilon)+\int\limits_{0}^{\kappa}\int\limits_{\mathcal{Z}}S_{\alpha}(\kappa-\upsilon)\Gamma_{4} (\upsilon, x_{\upsilon}, u) y(\upsilon)N(d\upsilon, du)\Big)d\upsilon\nonumber
	\end{align}
	\begin{align}
		&\quad+\int\limits_{0}^{\iota}S_{\alpha}(\iota-\upsilon)G\Big(\upsilon, x_{\upsilon}+y_{\upsilon}+ \int\limits_{0}^{\kappa}S_{\alpha}(\kappa-\upsilon)\Gamma_{1}(\upsilon, x_{\upsilon}) y(\upsilon)d\upsilon+ \int\limits_{0}^{\kappa}S_{\alpha}(\kappa-\upsilon)\Gamma_{2}(\upsilon, x_{\upsilon}) y(\upsilon)d\upsilon\nonumber\\
				&\quad+ \int\limits_{0}^{\kappa}S_{\alpha}(\kappa-\upsilon)\Gamma_{3}(\upsilon, x_{\upsilon}) y(\upsilon)dw(\upsilon)+\int\limits_{0}^{\kappa}\int\limits_{\mathcal{Z}}S_{\alpha}(\kappa-\upsilon)\Gamma_{4} (\upsilon, x_{\upsilon}, u) y(\upsilon)N(d\upsilon, du)\Big)dw(\upsilon)\nonumber\\
				&\quad+\int\limits_{0}^{\iota}\int\limits_{\mathcal{Z}}S_{\alpha}(\iota-\upsilon)\sigma\Big(\upsilon, x_{\upsilon}+y_{\upsilon}+ \int\limits_{0}^{\kappa}S_{\alpha}(\kappa-\upsilon)\Gamma_{1}(\upsilon, x_{\upsilon}) y(\upsilon)d\upsilon+ \int\limits_{0}^{\kappa}S_{\alpha}(\kappa-\upsilon)\Gamma_{2}(\upsilon, x_{\upsilon}) y(\upsilon)d\upsilon\nonumber\\
						&\quad+ \int\limits_{0}^{\kappa}S_{\alpha}(\kappa-\upsilon)\Gamma_{3}(\upsilon, x_{\upsilon}) y(\upsilon)dw(\upsilon)+\int\limits_{0}^{\kappa}\int\limits_{\mathcal{Z}}S_{\alpha}(\kappa-\upsilon)\Gamma_{4} (\upsilon, x_{\upsilon}, u) y(\upsilon)N(d\upsilon, du)\Big)N(d\upsilon, du).
\end{align}
$ \therefore  $ Using  \eqref{eq6} and H$ \ddot{o} $lder's inequality, we have
\begin{align}\label{eq7}
&\mathbb{E}\big\|y_{n+1}(\iota)\big\|^{p}\nonumber\\
&=\ 5^{p-1}\Big\{\mathbb{E}\Big\|-g(\iota,x_{\iota})+\int\limits_{0}^{\iota}S_{\alpha}(\iota-\upsilon)\Gamma_{1}(\upsilon, x_{\upsilon}) y(\upsilon)d\upsilon- \int\limits_{0}^{\iota}\mathcal{A}S_{\alpha}(\iota-\upsilon)g(\upsilon, x_{\upsilon})d\upsilon\nonumber\\
&\quad+\int\limits_{0}^{\iota}S_{\alpha}(\iota-\upsilon)\Gamma_{1}(\upsilon, x_{\upsilon}) y(\upsilon)d\upsilon+ \int\limits_{0}^{\iota}S_{\alpha}(\iota-\upsilon) f(\upsilon,x_{\upsilon})d\upsilon+\int\limits_{0}^{\iota}S_{\alpha}(\iota-\upsilon)\Gamma_{2}(\upsilon, x_{\upsilon}) y(\upsilon)d\upsilon \nonumber\\
&\quad+ \int\limits_{0}^{\iota}S_{\alpha}(\iota-\upsilon) G(\upsilon, x_{\upsilon})dw(\upsilon)+\int\limits_{0}^{\iota}S_{\alpha}(\iota-\upsilon)\Gamma_{2}(\upsilon, x_{\upsilon}) y(\upsilon)dw(\upsilon)\nonumber\\
&\quad+ \int\limits_{0}^{\iota}\int\limits_{\mathcal{Z}}S_{\alpha}(\iota-\upsilon) \sigma(\upsilon, x_{\upsilon}, u)N(d\upsilon, du)+\int\limits_{0}^{\iota}\int\limits_{\mathcal{Z}}S_{\alpha}(\iota-\upsilon)\Gamma_{4} (\upsilon, x_{\upsilon}, u) y(\upsilon)\lambda du\nonumber\\
	 &\quad-g\Big(\iota, x_{\iota}+ y_{\iota}+\int\limits_{0}^{\iota}S_{\alpha}(\iota-\upsilon)\Gamma_{1}(\upsilon, x_{\upsilon}) y(\upsilon)d\upsilon+\int\limits_{0}^{\iota}S_{\alpha}(\iota-\upsilon)\Gamma_{2}(\upsilon, x_{\upsilon}) y(\upsilon)d\upsilon\nonumber\\
	&\quad+ \int\limits_{0}^{\iota}S_{\alpha}(\iota-\upsilon)\Gamma_{3}(\upsilon, x_{\upsilon}) y(\upsilon)dw(\upsilon)+\int\limits_{0}^{\iota}\int\limits_{\mathcal{Z}}S_{\alpha}(\iota-\upsilon)\Gamma_{4} (\upsilon, x_{\upsilon}, u) y(\upsilon)\lambda du\Big)\nonumber\\
&\quad+	\int\limits_{0}^{\iota}\mathcal{A}S_{\alpha}(\iota-\upsilon)g\Big(\upsilon, x_{\upsilon}+y_{\upsilon}+ \int\limits_{0}^{\kappa}S_{\alpha}(\kappa-\upsilon)\Gamma_{1}(\upsilon, x_{\upsilon}) y(\upsilon)d\upsilon+ \int\limits_{0}^{\kappa}S_{\alpha}(\kappa-\upsilon)\Gamma_{2}(\upsilon, x_{\upsilon}) y(\upsilon)d\upsilon\nonumber\\
	&\quad+ \int\limits_{0}^{\kappa}S_{\alpha}(\kappa-\upsilon)\Gamma_{3}(\upsilon, x_{\upsilon}) y(\upsilon)dw(\upsilon)+\int\limits_{0}^{\kappa}\int\limits_{\mathcal{Z}}S_{\alpha}(\kappa-\upsilon)\Gamma_{4} (\upsilon, x_{\upsilon}, u) y(\upsilon)N(d\upsilon, du)\Big)d\upsilon\nonumber
\end{align}
\begin{align}
	&\quad+\int\limits_{0}^{\iota}S_{\alpha}(\iota-\upsilon)f\Big(\upsilon, x_{\upsilon}+y_{\upsilon}+ \int\limits_{0}^{\kappa}S_{\alpha}(\kappa-\upsilon)\Gamma_{1}(\upsilon, x_{\upsilon}) y(\upsilon)d\upsilon+ \int\limits_{0}^{\kappa}S_{\alpha}(\kappa-\upsilon)\Gamma_{2}(\upsilon, x_{\upsilon}) y(\upsilon)d\upsilon\nonumber\\
		&\quad+ \int\limits_{0}^{\kappa}S_{\alpha}(\kappa-\upsilon)\Gamma_{3}(\upsilon, x_{\upsilon}) y(\upsilon)dw(\upsilon)+\int\limits_{0}^{\kappa}\int\limits_{\mathcal{Z}}S_{\alpha}(\kappa-\upsilon)\Gamma_{4} (\upsilon, x_{\upsilon}, u) y(\upsilon)N(d\upsilon, du)\Big)d\upsilon\nonumber\\
		&\quad+\int\limits_{0}^{\iota}S_{\alpha}(\iota-\upsilon)G\Big(\upsilon, x_{\upsilon}+y_{\upsilon}+ \int\limits_{0}^{\kappa}S_{\alpha}(\kappa-\upsilon)\Gamma_{1}(\upsilon, x_{\upsilon}) y(\upsilon)d\upsilon+ \int\limits_{0}^{\kappa}S_{\alpha}(\kappa-\upsilon)\Gamma_{2}(\upsilon, x_{\upsilon}) y(\upsilon)d\upsilon\nonumber\\
				&\quad+ \int\limits_{0}^{\kappa}S_{\alpha}(\kappa-\upsilon)\Gamma_{3}(\upsilon, x_{\upsilon}) y(\upsilon)dw(\upsilon)+\int\limits_{0}^{\kappa}\int\limits_{\mathcal{Z}}S_{\alpha}(\kappa-\upsilon)\Gamma_{4} (\upsilon, x_{\upsilon}, u) y(\upsilon)N(d\upsilon, du)\Big)dw(\upsilon)\nonumber\\
				&\quad+\int\limits_{0}^{\iota}\int\limits_{\mathcal{Z}}S_{\alpha}(\iota-\upsilon)\sigma\Big(\upsilon, x_{\upsilon}+y_{\upsilon}+ \int\limits_{0}^{\kappa}S_{\alpha}(\kappa-\upsilon)\Gamma_{1}(\upsilon, x_{\upsilon}) y(\upsilon)d\upsilon+ \int\limits_{0}^{\kappa}S_{\alpha}(\kappa-\upsilon)\Gamma_{2}(\upsilon, x_{\upsilon}) y(\upsilon)d\upsilon\nonumber\\
						&\quad+ \int\limits_{0}^{\kappa}S_{\alpha}(\kappa-\upsilon)\Gamma_{3}(\upsilon, x_{\upsilon}) y(\upsilon)dw(\upsilon)+\int\limits_{0}^{\kappa}\int\limits_{\mathcal{Z}}S_{\alpha}(\kappa-\upsilon)\Gamma_{4} (\upsilon, x_{\upsilon}, u) y(\upsilon)N(d\upsilon, du)\Big)N(d\upsilon, du)\Big\|^{p}\nonumber\Big\}\\
&\leq \ 5^{p-1}\sum\limits_{i=1}^{5} a_{i}.
\end{align}
Now, we compute the R.H.S of \eqref{eq7}. By using  Definitions \ref{def2.9} (ii),
\begin{align*}
a_{1} &\ \leq \ \mathbb{E}\Big\|
	 \mathcal{A}^{\gamma}g\Big(\iota, x_{\iota}+ y_{\iota}+\int\limits_{0}^{\iota}S_{\alpha}(\iota-\upsilon)\Gamma_{1}(\upsilon, x_{\upsilon}) y(\upsilon)d\upsilon+\int\limits_{0}^{\iota}S_{\alpha}(\iota-\upsilon)\Gamma_{2}(\upsilon, x_{\upsilon}) y(\upsilon)d\upsilon\nonumber\\
	&\quad+ \int\limits_{0}^{\iota}S_{\alpha}(\iota-\upsilon)\Gamma_{3}(\upsilon, x_{\upsilon}) y(\upsilon)dw(\upsilon)+\int\limits_{0}^{\iota}\int\limits_{\mathcal{Z}}S_{\alpha}(\iota-\upsilon)\Gamma_{4} (\upsilon, x_{\upsilon}, u) y(\upsilon)\lambda du\Big)\\
	&\quad-\mathcal{A}^{\gamma}g(\iota, x_{\iota})- \int\limits_{0}^{\iota}S_{\alpha}(\iota-\upsilon)\Gamma_{1}(\upsilon, x_{\upsilon}) y(\upsilon)d\upsilon\Big\|^{p}\\
&\ \leq \ \|\mathcal{A}^{-\gamma}\|^{p} \ \hat{a}_{1}^{p}\ \mathbb{E}\|y_{\iota}\|^{p}. 
\end{align*}
By using Lemma \ref{2.6} and Definition \ref{def2.9} (i), we get the following estimation:
\begin{align*}
a_{2} &\ = \  \mathbb{E}\Big\|	\int\limits_{0}^{\iota}\mathcal{A}S_{\alpha}(\iota-\upsilon)g\Big(\upsilon, x_{\upsilon}+y_{\upsilon}+ \int\limits_{0}^{\kappa}S_{\alpha}(\kappa-\upsilon)\Gamma_{1}(\upsilon, x_{\upsilon}) y(\upsilon)d\upsilon+ \int\limits_{0}^{\kappa}S_{\alpha}(\kappa-\upsilon)\Gamma_{2}(\upsilon, x_{\upsilon}) y(\upsilon)d\upsilon\nonumber\\
	&\quad+ \int\limits_{0}^{\kappa}S_{\alpha}(\kappa-\upsilon)\Gamma_{3}(\upsilon, x_{\upsilon}) y(\upsilon)dw(\upsilon)+\int\limits_{0}^{\kappa}\int\limits_{\mathcal{Z}}S_{\alpha}(\kappa-\upsilon)\Gamma_{4} (\upsilon, x_{\upsilon}, u) y(\upsilon)N(d\upsilon, du)\Big)d\upsilon
\end{align*}
\begin{align*}
	&\quad-\int\limits_{0}^{\iota}\mathcal{A}S_{\alpha}(\iota-\upsilon)g(\upsilon, x_{\upsilon}) d\upsilon
	-\int\limits_{0}^{\iota}S_{\alpha}(\iota-\upsilon)\Gamma_{1}(\upsilon, x_{\upsilon})y(\upsilon) d\upsilon	\Big\|^{p}\\
	&\ \leq \	\mathbb{E}\Big\|	\int\limits_{0}^{\iota}\mathcal{A}^{1-\gamma}S_{\alpha}(\iota-\upsilon)\mathcal{A}^{\gamma}g\Big(\upsilon, x_{\upsilon}+y_{\upsilon}+ \int\limits_{0}^{\kappa}S_{\alpha}(\kappa-\upsilon)\Gamma_{1}(\upsilon, x_{\upsilon}) y(\upsilon)d\upsilon+ \int\limits_{0}^{\kappa}S_{\alpha}(\kappa-\upsilon)\Gamma_{2}(\upsilon, x_{\upsilon}) y(\upsilon)d\upsilon\nonumber\\
		&\quad+ \int\limits_{0}^{\kappa}S_{\alpha}(\kappa-\upsilon)\Gamma_{3}(\upsilon, x_{\upsilon}) y(\upsilon)dw(\upsilon)+\int\limits_{0}^{\kappa}\int\limits_{\mathcal{Z}}S_{\alpha}(\kappa-\upsilon)\Gamma_{4} (\upsilon, x_{\upsilon}, u) y(\upsilon)N(d\upsilon, du)\Big)d\upsilon\\
		&\quad-\int\limits_{0}^{\iota}\mathcal{A}^{1-\gamma}S_{\alpha}(\iota-\upsilon)\mathcal{A}^{\gamma}g(\upsilon, x_{\upsilon}) d\upsilon
		-\int\limits_{0}^{\iota}S_{\alpha}(\iota-\upsilon)\Gamma_{1}(\upsilon, x_{\upsilon})y(\upsilon) d\upsilon	\Big\|^{p}\\
	&\ \leq \ \alpha^{p}\ c_{\mu}^{p}\ \Big[\frac{\iota^{p\alpha \mu}}{p\alpha \mu}\Big]^{p-1} \int\limits_{0}^{\iota}(\iota-\upsilon)^{-p\alpha \mu} \hat{a}_{1}^{p} \mathbb{E}\|y_{n}(\upsilon)\|^{p}\\
	 	&\ \leq \ \alpha^{p}\ c_{\mu}^{p}\ \hat{a}_{1}^{p}\ \Big[\frac{\iota^{p\alpha \mu}}{p\alpha \mu}\Big]^{p-1} \int\limits_{0}^{\iota}(\iota-\upsilon)^{-p\alpha \mu}  \mathbb{E}\|y_{n}(\upsilon)\|^{p}.	
	\end{align*}
By using Lemma \ref{2.6} and Definition \ref{def2.9} (i), we get the following estimation
\begin{align*}
a_{3} &\ = \  \mathbb{E}\Big\|\int\limits_{0}^{\iota}S_{\alpha}(\iota-\upsilon)f\Big(\upsilon, x_{\upsilon}+y_{\upsilon}+ \int\limits_{0}^{\kappa}S_{\alpha}(\kappa-\upsilon)\Gamma_{1}(\upsilon, x_{\upsilon}) y(\upsilon)d\upsilon+ \int\limits_{0}^{\kappa}S_{\alpha}(\kappa-\upsilon)\Gamma_{2}(\upsilon, x_{\upsilon}) y(\upsilon)d\upsilon\nonumber\\
		&\quad+ \int\limits_{0}^{\kappa}S_{\alpha}(\kappa-\upsilon)\Gamma_{3}(\upsilon, x_{\upsilon}) y(\upsilon)dw(\upsilon)+\int\limits_{0}^{\kappa}\int\limits_{\mathcal{Z}}S_{\alpha}(\kappa-\upsilon)\Gamma_{4} (\upsilon, x_{\upsilon}, u) y(\upsilon)N(d\upsilon, du)\Big)d\upsilon\\
		&\quad-\int\limits_{0}^{\iota}S_{\alpha}(\iota-\upsilon)f(\upsilon, x_{\upsilon}) d\upsilon
			-\int\limits_{0}^{\iota}S_{\alpha}(\iota-\upsilon)\Gamma_{2}(\upsilon, x_{\upsilon})y(\upsilon) d\upsilon	\Big\|^{p}\\
	&	\leq\ \mathbb{E}\int\limits_{0}^{\iota}\Big\|S_{\alpha}(\iota-\upsilon)f\Big(\upsilon, x_{\upsilon}+y_{\upsilon}+ \int\limits_{0}^{\kappa}S_{\alpha}(\kappa-\upsilon)\Gamma_{1}(\upsilon, x_{\upsilon}) y(\upsilon)d\upsilon+ \int\limits_{0}^{\kappa}S_{\alpha}(\kappa-\upsilon)\Gamma_{2}(\upsilon, x_{\upsilon}) y(\upsilon)d\upsilon\nonumber\\
		&\quad+ \int\limits_{0}^{\kappa}S_{\alpha}(\kappa-\upsilon)\Gamma_{3}(\upsilon, x_{\upsilon}) y(\upsilon)dw(\upsilon)+\int\limits_{0}^{\kappa}\int\limits_{\mathcal{Z}}S_{\alpha}(\kappa-\upsilon)\Gamma_{4} (\upsilon, x_{\upsilon}, u) y(\upsilon)N(d\upsilon, du)\Big)d\upsilon\\
		&\quad-\int\limits_{0}^{\iota}S_{\alpha}(\iota-\upsilon)f(\upsilon, x_{\upsilon}) d\upsilon
					-\int\limits_{0}^{\iota}S_{\alpha}(\iota-\upsilon)\Gamma_{2}(\upsilon, x_{\upsilon})y(\upsilon) d\upsilon	\Big\|^{p}\\
		&	\leq\ M^{p} \ \iota^{p-1} \ b_{2}^{p}\ \int\limits_{0}^{\iota}(\iota-\upsilon)  \mathbb{E}\|y_{n}(\upsilon)\|^{p}.		
\end{align*}	
By using Lemma \ref{lm2} and Definitions \ref{2.6}, Definitions \ref{def2.9} (iii),  one can obtain	
\begin{align*}
a_{4} &\ = \ \mathbb{E}\Big\|\int\limits_{0}^{\iota}S_{\alpha}(\iota-\upsilon)G\Big(\upsilon, x_{\upsilon}+y_{\upsilon}+ \int\limits_{0}^{\kappa}S_{\alpha}(\kappa-\upsilon)\Gamma_{1}(\upsilon, x_{\upsilon}) y(\upsilon)d\upsilon+ \int\limits_{0}^{\kappa}S_{\alpha}(\kappa-\upsilon)\Gamma_{2}(\upsilon, x_{\upsilon}) y(\upsilon)d\upsilon\nonumber\\
&\quad+ \int\limits_{0}^{\kappa}S_{\alpha}(\kappa-\upsilon)\Gamma_{3}(\upsilon, x_{\upsilon}) y(\upsilon)dw(\upsilon)+\int\limits_{0}^{\kappa}\int\limits_{\mathcal{Z}}S_{\alpha}(\kappa-\upsilon)\Gamma_{4} (\upsilon, x_{\upsilon}, u) y(\upsilon)N(d\upsilon, du)\Big)dw(\upsilon)\\
&\quad-\int\limits_{0}^{\iota}S_{\alpha}(\iota-\upsilon)G(\upsilon, x_{\upsilon}) d\upsilon
-\int\limits_{0}^{\iota}S_{\alpha}(\iota-\upsilon)\Gamma_{3}(\upsilon, x_{\upsilon})y(\upsilon) dw(\upsilon)\Big\|^{p}\\
&\ \leq \ C_{p}\ M^{p} \Big[\int\limits_{0}^{\iota}\Big(\mathbb{E}\Big\|G\Big(\upsilon, x_{\upsilon}+y_{\upsilon}+ \int\limits_{0}^{\kappa}S_{\alpha}(\kappa-\upsilon)\Gamma_{1}(\upsilon, x_{\upsilon}) y(\upsilon)d\upsilon+ \int\limits_{0}^{\kappa}S_{\alpha}(\kappa-\upsilon)\Gamma_{2}(\upsilon, x_{\upsilon}) y(\upsilon)d\upsilon\nonumber\\
&\quad+ \int\limits_{0}^{\kappa}S_{\alpha}(\kappa-\upsilon)\Gamma_{3}(\upsilon, x_{\upsilon}) y(\upsilon)dw(\upsilon)+\int\limits_{0}^{\kappa}\int\limits_{\mathcal{Z}}S_{\alpha}(\kappa-\upsilon)\Gamma_{4} (\upsilon, x_{\upsilon}, u) y(\upsilon)N(d\upsilon, du)\Big)dw(\upsilon)\\
&\quad-\int\limits_{0}^{\iota}S_{\alpha}(\iota-\upsilon)G(\upsilon, x_{\upsilon}) d\upsilon
-\int\limits_{0}^{\iota}S_{\alpha}(\iota-\upsilon)\Gamma_{3}(\upsilon, x_{\upsilon})y(\upsilon) dw(\upsilon)\Big\|^{p}\Big)^{\frac{2}{p}}d\upsilon\Big]^{\frac{p}{2}}
&\ \leq \ C_{p}\  M^{P}\ \iota^{\frac{p}{2}-1} \ \hat{a}_{3}^{p}  \int\limits_{0}^{\iota}(\iota-\upsilon)  \mathbb{E}\|y_{n}(\upsilon)\|^{p}.	
	\end{align*}
By using Lemma \ref{2.8} and Definitions \ref{2.6}, Definitions \ref{def2.9} (iv),  one can obtain	
\begin{align*}
		a_{5} &\ = \ \int\limits_{0}^{\iota}\int\limits_{\mathcal{Z}}S_{\alpha}(\iota-\upsilon)\sigma\Big(\upsilon, x_{\upsilon}+y_{\upsilon}+ \int\limits_{0}^{\kappa}S_{\alpha}(\kappa-\upsilon)\Gamma_{1}(\upsilon, x_{\upsilon}) y(\upsilon)d\upsilon+ \int\limits_{0}^{\kappa}S_{\alpha}(\kappa-\upsilon)\Gamma_{2}(\upsilon, x_{\upsilon}) y(\upsilon)d\upsilon\nonumber\\
		&\quad+ \int\limits_{0}^{\kappa}S_{\alpha}(\kappa-\upsilon)\Gamma_{3}(\upsilon, x_{\upsilon}) y(\upsilon)dw(\upsilon)+\int\limits_{0}^{\kappa}\int\limits_{\mathcal{Z}}S_{\alpha}(\kappa-\upsilon)\Gamma_{4} (\upsilon, x_{\upsilon}, u) y(\upsilon)N(d\upsilon, du)\Big)N(d\upsilon, du)\\
		&\quad-\int\limits_{0}^{\iota}\int\limits_{\mathcal{Z}}S_{\alpha}(\iota-\upsilon)\sigma(\upsilon, x_{\upsilon}, u) N(du, d\upsilon)- \int\limits_{0}^{\iota}\int\limits_{\mathcal{Z}}S_{\alpha}(\iota-\upsilon)\Gamma_{4}(\upsilon, x_{\upsilon}, u) N(d\upsilon, du)\Big\|^{p}\\
		&\ \leq \ k(p)\ M^{p}\Big\{ \Big(\frac{\iota^{3}}{3}\Big)^{\frac{p}{2}}\hat{a}_{4}^{p}+ \Big(\frac{\iota^{2p+1}}{2p+1}\Big)^{\frac{1}{2}}\hat{a}_{5}^{p}\Big\} \mathbb{E}\|y_{n}(\upsilon)\|^{p}.	
	\end{align*}	
Combining these results together with equation \eqref {eq7},
\begin{align}\label{eq9}
	\mathbb{E}\|y_{n+1}(\iota)\| &\leq \ 5^{p-1} \Big\{\|\mathcal{A}^{-\gamma}\|^{p} \ \hat{a}_{1}^{p}+\alpha^{p}\ c_{\mu}^{p}\ \hat{a}_{1}^{p}\ \Big[\frac{\iota^{p\alpha \mu}}{p\alpha \mu}\Big]^{p}+k(p)\ M^{p}\Big\{ \Big(\frac{\iota^{3}}{3}\Big)^{\frac{p}{2}}\hat{a}_{4}^{p}+ \Big(\frac{\iota^{2p+1}}{2p+1}\Big)^{\frac{1}{2}}\hat{a}_{5}^{p}\Big\}\Big\} \mathbb{E}\|y(\iota)\|^{p} \nonumber\\
	&\quad+5^{p-1}\Big\{M^{p} \ \iota^{p-1} \  \hat{a}_{2}^{p}+ C_{p}\  M^{P}\ \iota^{\frac{p}{2}-1} \ \hat{a}_{3}^{p}\Big\}\int\limits_{0}^{\iota}(\iota-\upsilon)  \mathbb{E}\|y_{n}(\upsilon)\|^{p}.
\end{align}
Now, we have to substitute the value in $ y_{n}(\iota) $ in above equation \eqref{eq9} to obtain
\begin{align*}
\mathbb{E}\|y_{n+1}(\iota)\| &\leq \ \Big(5^{p-1} \Big\{\|\mathcal{A}^{-\gamma}\|^{p} \ \hat{a}_{1}^{p}+\alpha^{p}\ c_{\mu}^{p}\ \hat{a}_{1}^{p}\ \Big[\frac{\iota^{p\alpha \mu}}{p\alpha \mu}\Big]^{p}+k(p)\ M^{p}\Big\{ \Big(\frac{\iota^{3}}{3}\Big)^{\frac{p}{2}}\hat{a}_{4}^{p}+ \Big(\frac{\iota^{2p+1}}{2p+1}\Big)^{\frac{1}{2}}\hat{a}_{5}^{p}\Big\}\Big\}\Big)^{2}\mathbb{E}\|y(\iota)\|^{p}\\
&\quad+\Big(5^{p-1} \Big\{M^{p} \ \iota^{p-1} \ \hat{a}_{2}^{p}+ C_{p}\  M^{P}\ \iota^{\frac{p}{2}-1} \ \hat{a}_{3}^{p}\Big\}\Big)^{2}
\int\limits_{0}^{\upsilon}\int\limits_{0}^{\iota}(\iota-\upsilon)  \mathbb{E}\|y_{n-1}(\upsilon_{1})\|^{p}d\upsilon_{1}d\upsilon.\\
&\leq \ \Big(5^{p-1} \Big\{\|\mathcal{A}^{-\gamma}\|^{p} \ \hat{a}_{1}^{p}+\alpha^{p}\ c_{\mu}^{p}\ \hat{a}_{1}^{p}\ \Big[\frac{\iota^{p\alpha \mu}}{p\alpha \mu}\Big]^{p}+k(p)\ M^{p}\Big\{ \Big(\frac{\iota^{3}}{3}\Big)^{\frac{p}{2}}\hat{a}_{4}^{p}+ \Big(\frac{\iota^{2p+1}}{2p+1}\Big)^{\frac{1}{2}}\hat{a}_{5}^{p}\Big\}\Big\}\Big)^{3}\mathbb{E}\|y(\iota)\|^{p}\\
&\quad+\Big(5^{p-1} \Big\{M^{p} \ \iota^{p-1} \ \hat{a}_{2}^{p}+ C_{p}\  M^{P}\ \iota^{\frac{p}{2}-1} \ \hat{a}_{3}^{p}\Big\}\Big)^{3}
 \int\limits_{0}^{\iota}\int\limits_{0}^{\upsilon}\int\limits_{0}^{\upsilon_{1}}(\iota-\upsilon)^{\alpha-1}  \mathbb{E}\|y_{n-2}(\upsilon_{2})\|^{p} d\upsilon_{2}d\upsilon_{1}d\upsilon.
 \end{align*}
From the above procedure and applying the mathematical induction method, we get
\begin{align*}
\mathbb{E}\|y_{n+1}(\iota)\| &\leq \ \Big(5^{p-1} \Big\{\|\mathcal{A}^{-\gamma}\|^{p} \ \hat{a}_{1}^{p}+\alpha^{p}\ c_{\mu}^{p}\ \hat{a}_{1}^{p}\ \Big[\frac{\iota^{p\alpha \mu}}{p\alpha \mu}\Big]^{p}+k(p)\ M^{p}\Big\{ \Big(\frac{\iota^{3}}{3}\Big)^{\frac{p}{2}}\hat{a}_{4}^{p}+ \Big(\frac{\iota^{2p+1}}{2p+1}\Big)^{\frac{1}{2}}\hat{a}_{5}^{p}\Big\}\Big\}\Big)^{n+1}\\
&\quad\times\mathbb{E}\|y_{0}(\upsilon)\|^{p}+ \Big(5^{p-1} \Big\{M^{p} \ \iota^{p-1} \ \hat{a}_{2}^{p}+ C_{p}\  M^{P}\ \iota^{\frac{p}{2}-1} \ \hat{a}_{3}^{p}\Big\}\Big)^{n+1}
 \frac{1}{n!}\int\limits_{0}^{\iota}(\iota-\upsilon)^{n}  \mathbb{E}\|y_{0}(\upsilon)\|^{p}d\upsilon.\\
&\leq \ \Big(5^{p-1} \Big\{\|\mathcal{A}^{-\gamma}\|^{p} \ \hat{a}_{1}^{p}+\alpha^{p}\ c_{\mu}^{p}\ \hat{a}_{1}^{p}\ \Big[\frac{\iota^{p\alpha \mu}}{p\alpha \mu}\Big]^{p}+k(p)\ M^{p}\Big\{ \Big(\frac{\iota^{3}}{3}\Big)^{\frac{p}{2}}\hat{a}_{4}^{p}+ \Big(\frac{\iota^{2p+1}}{2p+1}\Big)^{\frac{1}{2}}\hat{a}_{5}^{p}\Big\}\Big\}\Big)^{n+1}\\
&\quad\times\mathbb{E}\|y_{0}(\upsilon)\|^{p}+ \Big(5^{p-1} \Big\{M^{p} \ \iota^{p-1} \ \hat{a}_{2}^{p}+ C_{p}\  M^{P}\ \iota^{\frac{p}{2}-1} \ \hat{a}_{3}^{p}\Big\}\Big)^{n+1}
 \frac{1}{n!} \ \frac{\iota^{n+1}}{n+1} \mathbb{E}\|y_{0}(\upsilon)\|^{p};
\end{align*}
as $n\rightarrow \infty , \quad y_{n}(\iota) $  converges to zero in $ \mathbb{H}.$

\noindent
Case (ii): 
Now, we have to prove that $\{ x_{n}\}_{n=1}^{\infty}$ converges to the given solution of Cauchy problem \eqref{eq1}. \\
Now, we have to constructed the following:
\begin{align*}
x_{n+1}(\iota) &\ =\ x_{n}(\iota)-y_{n}(\iota)- \int\limits_{0}^{\iota}S_{\alpha}(\iota-\upsilon)\Gamma_{1}(\upsilon, x^{n}_{\upsilon}) y^{n}(\upsilon)d\upsilon- \int\limits_{0}^{\iota}S_{\alpha}(\iota-\upsilon)\Gamma_{2}(\upsilon, x^{n}_{\upsilon}) y^{n}(\upsilon)d\upsilon\\
&\quad-\int\limits_{0}^{\iota}S_{\alpha}(\iota-\upsilon)\Gamma_{3}(\upsilon, x^{n}_{\upsilon}) y^{n}(\upsilon)dw(\upsilon)- \int\limits_{0}^{\iota} \int\limits_{\mathcal{Z}}S_{\alpha}(\iota-\upsilon)\Gamma_{4}(\upsilon, x^{n}_{\upsilon}, u) y^{n}(\upsilon)N(d\upsilon, du).\\
x_{n+1}(\iota)- x_{n}(\iota)&\ =\ -y_{n}(\iota)- \int\limits_{0}^{\iota}S_{\alpha}(\iota-\upsilon)\Gamma_{1}(\upsilon, x^{n}_{\upsilon}) y^{n}(\upsilon)d\upsilon- \int\limits_{0}^{\iota}S_{\alpha}(\iota-\upsilon)\Gamma_{2}(\upsilon, x^{n}_{\upsilon}) y^{n}(\upsilon)d\upsilon\\
&\quad-\int\limits_{0}^{\iota}S_{\alpha}(\iota-\upsilon)\Gamma_{3}(\upsilon, x^{n}_{\upsilon}) y^{n}(\upsilon)dw(\upsilon)- \int\limits_{0}^{\iota} \int\limits_{\mathcal{Z}}S_{\alpha}(\iota-\upsilon)\Gamma_{4}(\upsilon, x^{n}_{\upsilon}, u) y^{n}(\upsilon)N(d\upsilon, du).
\end{align*}
\begin{align*}
\mathbb{E}\|x_{n+1}(\iota)- x_{n}(\iota)\|^{p}&\ =\ \mathbb{E}\Big\|-y_{n}(\iota)- \int\limits_{0}^{\iota}S_{\alpha}(\iota-\upsilon)\Gamma_{1}(\upsilon, x^{n}_{\upsilon}) y^{n}(\upsilon)d\upsilon- \int\limits_{0}^{\iota}S_{\alpha}(\iota-\upsilon)\Gamma_{2}(\upsilon, x^{n}_{\upsilon}) y^{n}(\upsilon)d\upsilon\\
&\quad-\int\limits_{0}^{\iota}S_{\alpha}(\iota-\upsilon)\Gamma_{3}(\upsilon, x^{n}_{\upsilon}) y^{n}(\upsilon)dw(\upsilon)- \int\limits_{0}^{\iota} \int\limits_{\mathcal{Z}}S_{\alpha}(\iota-\upsilon)\Gamma_{4}(\upsilon, x^{n}_{\upsilon}, u) y^{n}(\upsilon)N(d\upsilon, du)\Big\|^{p}\\
&\ \leq \ \frac{5^{p-1}}{(n-1)!} \Big[\Big(5^{p-1} \Big\{\|\mathcal{A}^{-\gamma}\|^{p} \ \hat{a}_{1}^{p}+\alpha^{p}\ c_{\mu}^{p}\ \hat{a}_{1}^{p}\ \Big[\frac{\iota^{p\alpha \mu}}{p\alpha \mu}\Big]^{p}\\
&\quad+k(p)\ M^{p}\Big\{ \Big(\frac{\iota^{3}}{3}\Big)^{\frac{p}{2}}\hat{a}_{4}^{p}+ \Big(\frac{\iota^{2p+1}}{2p+1}\Big)^{\frac{1}{2}}\hat{a}_{5}^{p}\Big\}\Big\}\Big)^{n}\\
&\quad+\Big(5^{p-1} \Big\{M^{p} \ \iota^{p-1} \ \hat{a}_{2}^{p}+ C_{p}\  M^{P}\ \iota^{\frac{p}{2}-1} \ \hat{a}_{3}^{p}\Big\}\Big)^{n}\frac{\iota^{n+1}}{n+1}\\
&\quad \times\Big(1+M^{p} \iota^{p}\ c_{1}^{p}+ C_{p} \ M^{p} \iota^{\frac{p}{2}}\ c_{2}^{p}+M^{p} k(p)[\iota^{\frac{p}{2}}c_{4}^{p}+\iota^{p}\hat{c}_{4}^{p}] \Big)\mathbb{E}\|y_{0}(\iota)\|^{p}\Big],
\end{align*}
for a given $m>0 $ s.t $ m < n,$ then from the above procedure, we have
\begin{align}\label{eq10}
\mathbb{E}\|x_{n}(\iota)- x_{m}(\iota)\|^{p}&\ \leq \ \mathbb{E}\|x_{n}(\iota)- x_{n-1}(\iota)\|^{p}+ \mathbb{E}\|x_{n-1}(\iota)- x_{n-2}(\iota)\|^{p}+\cdots+ \mathbb{E}\|x_{m+1}(\iota)- x_{m}(\iota)\|^{p}\nonumber\\
&\ \leq \ 5^{p-1} \Big(1+M^{p} \iota^{p}\ c_{1}^{p}+ C_{p} \ M^{p} \iota^{\frac{p}{2}}\ c_{2}^{p}+M^{p} k(p)[\iota^{\frac{p}{2}}c_{4}^{p}+\iota^{p}\hat{c}_{4}^{p}] \Big)\mathbb{E}\|y_{0}(\iota)\|^{p}\nonumber\\
&\quad \times\sum\limits_{k=m}^{n-1} \frac{1}{(k-1)!} \ A^{k},
\end{align}
where 
\begin{align*}
A^{k}&:= \ \frac{5^{p-1}}{(n-1)!} \Big[\Big(5^{p-1} \Big\{\|\mathcal{A}^{-\gamma}\|^{p} \ \hat{a}_{1}^{p}+\alpha^{p}\ c_{\mu}^{p}\ \hat{a}_{1}^{p}\ \Big[\frac{\iota^{p\alpha \mu}}{p\alpha \mu}\Big]^{p}
+k(p)\ M^{p}\Big\{ \Big(\frac{\iota^{3}}{3}\Big)^{\frac{p}{2}}\hat{a}_{4}^{p}+ \Big(\frac{\iota^{2p+1}}{2p+1}\Big)^{\frac{1}{2}}\hat{a}_{5}^{p}\Big\}\Big\}\Big)^{n}\\
&\quad+\Big(5^{p-1} \Big\{M^{p} \ \iota^{p-1} \ \hat{a}_{2}^{p}+ C_{p}\  M^{P}\ \iota^{\frac{p}{2}-1} \ \hat{a}_{3}^{p}\Big\}\Big)^{n}\frac{\iota^{n+1}}{n+1}.
 \end{align*}

\noindent 
Clearly, the R.H.S of \eqref{eq10} is a convergent series. Thus, $ \{x_{n}\}_{n=1}^{\infty} $ is a Cauchy sequence in $ \mathcal{B} $ and hence convergent uniformly to $ x^{*}$. Hence, 
\begin{align*}
y_{n}(\iota)&\ = \ x_{n}(\iota)-g(\iota, x_{\iota}^{n})- \int\limits_{0}^{\iota}\mathcal{A}S_{\alpha}(\iota-\upsilon)g(\upsilon, x_{\upsilon}^{n})d\upsilon+ \int\limits_{0}^{\iota}S_{\alpha}(\iota-\upsilon) f(\upsilon, x_{\upsilon}^{n})d\upsilon\\
&\quad+ \int\limits_{0}^{\iota}S_{\alpha}(\iota-\upsilon) G(\upsilon, x_{\upsilon}^{n})dw(\upsilon) 
	+ \int\limits_{0}^{\iota}\int\limits_{\mathcal{Z}}S_{\alpha}(\iota-\upsilon) \sigma(\upsilon, x_{\upsilon}^{n}, u)N(d\upsilon, du)-x_{0}(\iota).\\
\lim\limits_{n \rightarrow \infty}y_{n}(\iota)&\ = \ \lim\limits_{n \rightarrow \infty}x_{n}(\iota)- x_{0}(\iota)-\lim\limits_{n \rightarrow \infty}g(\iota, x_{\iota}^{n})-\lim\limits_{n \rightarrow \infty} \int\limits_{0}^{\iota}\mathcal{A}S_{\alpha}(\iota-\upsilon)g(\upsilon, x_{\upsilon}^{n})d\upsilon\\
&\quad+\lim\limits_{n \rightarrow \infty} \int\limits_{0}^{\iota}S_{\alpha}(\iota-\upsilon) f(\upsilon, x_{\upsilon}^{n})d\upsilon 
	+ \lim\limits_{n \rightarrow \infty}\int\limits_{0}^{\iota}S_{\alpha}(\iota-\upsilon) G(\upsilon, x_{\upsilon}^{n})dw(\upsilon)
\end{align*}
\begin{align*}
	&\quad+ \lim\limits_{n \rightarrow \infty}\int\limits_{0}^{\iota}\int\limits_{\mathcal{Z}}S_{\alpha}(\iota-\upsilon) \sigma(\upsilon, x_{\upsilon}^{n}, u)N(d\upsilon, du).\\
x^{*}(\iota) &\ = x_{0}(\iota)-g(\iota, x_{\iota}^{*})- \int\limits_{0}^{\iota}\mathcal{A}S_{\alpha}(\iota-\upsilon)g(\upsilon, x_{\upsilon}^{*})d\upsilon+ \int\limits_{0}^{\iota}S_{\alpha}(\iota-\upsilon) f(\upsilon, x_{\upsilon}^{*})d\upsilon \nonumber\\
	&\quad
	+ \int\limits_{0}^{\iota}S_{\alpha}(\iota-\upsilon) G(\upsilon, x_{\upsilon}^{*})dw(\upsilon)+ \int\limits_{0}^{\iota}\int\limits_{\mathcal{Z}}S_{\alpha}(\iota-\upsilon) \sigma(\upsilon, x_{\upsilon}^{*}, u)N(d\upsilon, du).
	\end{align*}
	\begin{align*}
x^{*}(\iota) &\ = C_{\alpha}(\iota)[\phi+g(0, \phi)] +S_{\alpha}(\iota)\eta -g(\iota, x_{\iota})- \int\limits_{0}^{\iota}\mathcal{A}S_{\alpha}(\iota-\upsilon)g(\upsilon, x_{\upsilon})d\upsilon+ \int\limits_{0}^{\iota}S_{\alpha}(\iota-\upsilon) f(\upsilon, x_{\upsilon})d\upsilon\nonumber\\
&\quad 
+ \int\limits_{0}^{\iota}S_{\alpha}(\iota-\upsilon) G(\upsilon, x_{\upsilon})dw(\upsilon)
+ \int\limits_{0}^{\iota}\int\limits_{\mathcal{Z}}S_{\alpha}(\iota-\upsilon) \sigma(\upsilon, x_{\upsilon}, u)N(d\upsilon, du).
\end{align*}
Hence, $ x^{*}(\iota) $ is the mild solution of the given Cauchy problem \eqref{eq1}.\\
\textbf{Uniqueness Result for Mild Solution:}\\
Next, our aim is to prove the uniqueness of the solution by utilizing the regularity property of the integral contractor. Let $ x_{1} $ and $ x_{2} $ be the mild solutions of the given Cauchy problem \eqref{eq1}. Then, we have\\
\noindent By using the regularity condition in Definition \ref{def2.11} with $ x=x_{1} $ and $ A= x_{2}-x_{1} $ s.t 
\begin{align*}
&y(\iota)+\int\limits_{0}^{\iota}S_{\alpha}(\iota-\upsilon)\Gamma_{1} (\upsilon, x(\upsilon))y(\upsilon) +\int\limits_{0}^{\iota}S_{\alpha}(\iota-\upsilon)\Gamma_{2} (\upsilon, x(\upsilon))y(\upsilon)\nonumber\\
&+\int\limits_{0}^{\iota}S_{\alpha}(\iota-\upsilon)\Gamma_{3} (\upsilon, x(\upsilon))y(\upsilon)dw(\upsilon)+\int\limits_{0}^{\iota}\int\limits_{\mathcal{Z}}S_{\alpha}(\iota-\upsilon)\Gamma_{4} (\upsilon, x(\upsilon), u)y(\upsilon)\lambda du= A(\iota)\\
&y(\iota)+\int\limits_{0}^{\iota}S_{\alpha}(\iota-\upsilon)\Gamma_{1} (\upsilon, x(\upsilon))y(\upsilon) +\int\limits_{0}^{\iota}S_{\alpha}(\iota-\upsilon)\Gamma_{2} (\upsilon, x(\upsilon))y(\upsilon)\nonumber\\
&+\int\limits_{0}^{\iota}S_{\alpha}(\iota-\upsilon)\Gamma_{3} (\upsilon, x(\upsilon))y(\upsilon)dw(\upsilon)+\int\limits_{0}^{\iota}\int\limits_{\mathcal{Z}}S_{\alpha}(\iota-\upsilon)\Gamma_{4} (\upsilon, x(\upsilon), u)y(\upsilon)N(du, d\upsilon)   
= x_{2} (\iota)-x_{1} (\iota).
\end{align*}
Hence, 
\begin{align}\label{eq11}
  x_{2}(\iota) &\ = \ x_{1}(\iota)+ y(\iota)+\int\limits_{0}^{\iota}S_{\alpha}(\iota-\upsilon)\Gamma_{1} (\upsilon, x(\upsilon))y(\upsilon)d\upsilon +\int\limits_{0}^{\iota}S_{\alpha}(\iota-\upsilon)\Gamma_{2} (\upsilon, x(\upsilon))y(\upsilon)d\upsilon\nonumber\\
  &+\int\limits_{0}^{\iota}S_{\alpha}(\iota-\upsilon)\Gamma_{3} (\upsilon, x(\upsilon))y(\upsilon)dw(\upsilon)+\int\limits_{0}^{\iota}\int\limits_{\mathcal{Z}}S_{\alpha}(\iota-\upsilon)\Gamma_{4} (\upsilon, x(\upsilon), u)y(\upsilon)N(du, d\upsilon).  
\end{align}
Since $ x_{1} $ , $ x_{2} $ are the mild solutions of the fractional  neutral stochastic system \eqref{eq1}, we have
\begin{align}\label{eq12}
x_{1}(\iota) &\ =\ -g(\iota, x_{\iota}^{1})- \int\limits_{0}^{\iota}\mathcal{A}S_{\alpha}(\iota-\upsilon)g(\upsilon, x_{\upsilon}^{1})d\upsilon+ \int\limits_{0}^{\iota}S_{\alpha}(\iota-\upsilon) f(\upsilon, x_{\upsilon}^{1})d\upsilon\nonumber\\
&\quad 
+ \int\limits_{0}^{\iota}S_{\alpha}(\iota-\upsilon) G(\upsilon, x_{\upsilon}^{1})dw(\upsilon)
+ \int\limits_{0}^{\iota}\int\limits_{\mathcal{Z}}S_{\alpha}(\iota-\upsilon) \sigma(\upsilon, x_{\upsilon}^{1}, u)N(d\upsilon, du).\nonumber\\ 
x_{2}(\iota) &\ =\ -g(\iota, x_{\iota}^{2})- \int\limits_{0}^{\iota}\mathcal{A}S_{\alpha}(\iota-\upsilon)g(\upsilon, x_{\upsilon}^{2})d\upsilon+ \int\limits_{0}^{\iota}S_{\alpha}(\iota-\upsilon) f(\upsilon, x_{\upsilon}^{2})d\upsilon\nonumber\\
&\quad 
+ \int\limits_{0}^{\iota}S_{\alpha}(\iota-\upsilon) G(\upsilon, x_{\upsilon}^{2})dw(\upsilon)
+ \int\limits_{0}^{\iota}\int\limits_{\mathcal{Z}}S_{\alpha}(\iota-\upsilon) \sigma(\upsilon, x_{\upsilon}^{2}, u)N(d\upsilon, du).\nonumber\\
 x_{2}(\iota)-x_{1}(\iota)&\ =\ g(\iota, x_{\iota}^{2})-g(\iota, x_{\iota}^{1})+\int\limits_{0}^{\iota}\mathcal{A}S_{\alpha}(\iota-\upsilon)\big[g(\upsilon, x_{\upsilon}^{2})-g(\upsilon, x_{\upsilon}^{1})\big]d\upsilon\nonumber\\
&\quad+\int\limits_{0}^{\iota}S_{\alpha}(\iota-\upsilon)\big[f(\upsilon, x_{\upsilon}^{2})-f(\upsilon, x_{\upsilon}^{1})\big]d\upsilon+ \int\limits_{0}^{\iota} S_{\alpha}(\iota-\upsilon) \Big[G(\upsilon, x_{\upsilon}^{2})- G(\upsilon, x_{\upsilon}^{1})\Big] dw(\upsilon)\nonumber\\
&\quad + \int\limits_{0}^{\iota}\int\limits_{\mathcal{Z}} S_{\alpha}(\iota-\upsilon)\Big[\sigma(\upsilon, x_{\upsilon}^{2},u)-\sigma(\upsilon, x_{\upsilon}^{1}, u)\Big]N(d\upsilon, du).
\end{align}
Now we have to substitute equation \eqref{eq11} in \eqref{eq12}, we get
\begin{align}
&\mathbb{E}\|x_{2}(\iota)-x_{1}(\iota)\|^{p}\nonumber\\
&\ =\mathbb{E}\Big\|g\Big(\iota, x_{1}(\iota)+ y(\iota)+\int\limits_{0}^{\iota}S_{\alpha}(\iota-\upsilon)\Gamma_{1} (\upsilon, x(\upsilon))y(\upsilon)d\upsilon +\int\limits_{0}^{\iota}S_{\alpha}(\iota-\upsilon)\Gamma_{2} (\upsilon, x(\upsilon))y(\upsilon)d\upsilon\nonumber\\
  &\quad+\int\limits_{0}^{\iota}S_{\alpha}(\iota-\upsilon)\Gamma_{3} (\upsilon, x(\upsilon))y(\upsilon)dw(\upsilon)+\int\limits_{0}^{\iota}\int\limits_{\mathcal{Z}}S_{\alpha}(\iota-\upsilon)\Gamma_{4} (\upsilon, x(\upsilon), u)y(\upsilon)N(du, d\upsilon)\Big)-g(\iota, x_{\iota}^{1})\nonumber
\end{align}
\begin{align}
    &\quad-\int\limits_{0}^{\iota}S_{\alpha}(\iota-\upsilon)\Gamma_{1}(\upsilon, x_{1}(\upsilon)) y(\upsilon)d\upsilon+\int\limits_{0}^{\iota}\mathcal{A}S_{\alpha}(\iota-\upsilon)\Big[g(\upsilon,  x_{1}(\upsilon)+ y(\upsilon)+\int\limits_{0}^{\kappa}S_{\alpha}(\kappa-\upsilon)\Gamma_{1} (\upsilon, x(\upsilon))y(\upsilon)d\upsilon\nonumber\\
    &\quad +\int\limits_{0}^{\kappa}S_{\alpha}(\kappa-\upsilon)\Gamma_{2} (\upsilon, x(\upsilon))y(\upsilon)d\upsilon+\int\limits_{0}^{\kappa}S_{\alpha}(\kappa-\upsilon)\Gamma_{3} (\upsilon, x(\upsilon))y(\upsilon)dw(\upsilon)\nonumber
      &\quad+\int\limits_{0}^{\kappa}\int\limits_{\mathcal{Z}}S_{\alpha}(\kappa-\upsilon)\Gamma_{4} (\upsilon, x(\upsilon), u)y(\upsilon)N(du, d\upsilon))(\upsilon)-g(\upsilon, x_{\upsilon}^{1})-\int\limits_{0}^{\iota}S_{\alpha}(\iota-\upsilon)\Gamma_{1} (\upsilon, x(\upsilon))y(\upsilon)d\upsilon\Big]d\upsilon\nonumber\\
      &\quad+\int\limits_{0}^{\iota}S_{\alpha}(\iota-\upsilon)\Big[f(\upsilon,  x_{1}(\upsilon)+ y(\upsilon)+\int\limits_{0}^{\kappa}S_{\alpha}(\kappa-\upsilon)\Gamma_{1} (\upsilon, x(\upsilon))y(\upsilon)d\upsilon\nonumber\\
      &\quad +\int\limits_{0}^{\kappa}S_{\alpha}(\kappa-\upsilon)\Gamma_{2} (\upsilon, x(\upsilon))y(\upsilon)d\upsilon+\int\limits_{0}^{\kappa}S_{\alpha}(\kappa-\upsilon)\Gamma_{3} (\upsilon, x(\upsilon))y(\upsilon)dw(\upsilon)\nonumber\\
          &\quad+\int\limits_{0}^{\kappa}\int\limits_{\mathcal{Z}}S_{\alpha}(\kappa-\upsilon)\Gamma_{4} (\upsilon, x(\upsilon), u)y(\upsilon)N(du, d\upsilon))(\upsilon)-f(\upsilon, x_{\upsilon}^{1})-\int\limits_{0}^{\iota}S_{\alpha}(\iota-\upsilon)\Gamma_{2} (\upsilon, x(\upsilon))y(\upsilon)d\upsilon\Big]d\upsilon\nonumber\\
          &\quad+\int\limits_{0}^{\iota}S_{\alpha}(\iota-\upsilon)\Big[G(\upsilon,  x_{1}(\upsilon)+ y(\upsilon)+\int\limits_{0}^{\kappa}S_{\alpha}(\kappa-\upsilon)\Gamma_{1} (\upsilon, x(\upsilon))y(\upsilon)d\upsilon\nonumber\\
         &\quad +\int\limits_{0}^{\kappa}S_{\alpha}(\kappa-\upsilon)\Gamma_{2} (\upsilon, x(\upsilon))y(\upsilon)d\upsilon+\int\limits_{0}^{\kappa}S_{\alpha}(\kappa-\upsilon)\Gamma_{3} (\upsilon, x(\upsilon))y(\upsilon)dw(\upsilon)\nonumber\\
         &\quad+\int\limits_{0}^{\kappa}\int\limits_{\mathcal{Z}}S_{\alpha}(\kappa-\upsilon)\Gamma_{4} (\upsilon, x(\upsilon), u)y(\upsilon)N(du, d\upsilon))(\upsilon)-G(\upsilon, x_{\upsilon}^{1})-\int\limits_{0}^{\iota}S_{\alpha}(\iota-\upsilon)\Gamma_{3} (\upsilon, x(\upsilon))y(\upsilon)dw(\upsilon)\Big]dw(\upsilon)\nonumber\\
     &\quad+\int\limits_{0}^{\iota}\int\limits_{\mathcal{Z}}S_{\alpha}(\iota-\upsilon)\Big[\sigma(\upsilon,  x_{1}(\upsilon)+ y(\upsilon)+\int\limits_{0}^{\kappa}S_{\alpha}(\kappa-\upsilon)\Gamma_{1} (\upsilon, x(\upsilon))y(\upsilon)d\upsilon\nonumber\\
     &\quad +\int\limits_{0}^{\kappa}S_{\alpha}(\kappa-\upsilon)\Gamma_{2} (\upsilon,x(\upsilon))y(\upsilon)d\upsilon+\int\limits_{0}^{\kappa}S_{\alpha}(\kappa-\upsilon)\Gamma_{3} (\upsilon, x(\upsilon))y(\upsilon)dw(\upsilon)\nonumber\\
         &\quad+\int\limits_{0}^{\kappa}\int\limits_{\mathcal{Z}}S_{\alpha}(\kappa-\upsilon)\Gamma_{4} (\upsilon, x(\upsilon), u)y(\upsilon)N(du, d\upsilon), u)(\upsilon)-\sigma(\upsilon, x_{\upsilon}^{1}, u)\nonumber\\
        &\quad-\int\limits_{0}^{\iota}\int\limits_{\mathcal{Z}}S_{\alpha}(\iota-\upsilon)\Gamma_{4} (\upsilon, x(\upsilon), u)y(\upsilon)N(d\upsilon, du)\Big]N(d\upsilon, du).  
          \end{align}
 Hence, 
\begin{align*}
\mathbb{E}\|A(\iota)\|^{p} & \ \leq \ 5^{p-1} \Big\{\|\mathcal{A}^{-\gamma}\|^{p} \ \hat{a}_{1}^{p}+\alpha^{p}\ c_{\mu}^{p}\ \hat{a}_{1}^{p}\ \Big[\frac{\iota^{p\alpha \mu}}{p\alpha \mu}\Big]^{p}+k(p)\ M^{p}\\
&\quad\times\Big\{ \Big(\frac{\iota^{3}}{3}\Big)^{\frac{p}{2}}\hat{a}_{4}^{p}+ \Big(\frac{\iota^{2p+1}}{2p+1}\Big)^{\frac{1}{2}}\hat{a}_{5}^{p}\Big\}\Big\} \mathbb{E}\|A(\iota)\|^{p} \nonumber
\end{align*}
\begin{align*}
	&\quad+5^{p-1}\Big\{M^{p} \ \iota^{p-1} \  \hat{a}_{2}^{p}+ C_{p}\  M^{P}\ \iota^{\frac{p}{2}-1} \ \hat{a}_{3}^{p}\Big\}\int\limits_{0}^{\iota}(\iota-\upsilon)  \mathbb{E}\|A(\upsilon)\|^{p}.\nonumber\\
	& :=\ \oplus_{1}\mathbb{E}\|A(\iota)\|^{p}+\oplus_{2}\int\limits_{0}^{\iota}
\mathbb{E}\|A(\upsilon)\|^{p}d\upsilon,
\end{align*}
where
\begin{align*}
&  \oplus_{1}:=5^{p-1} \Big\{\|\mathcal{A}^{-\gamma}\|^{p} \ \hat{a}_{1}^{p}+\alpha^{p}\ c_{\mu}^{p}\ \hat{a}_{1}^{p}\ \Big[\frac{\iota^{p\alpha \mu}}{p\alpha \mu}\Big]^{p}+k(p)\ M^{p}\Big\{ \Big(\frac{\iota^{3}}{3}\Big)^{\frac{p}{2}}\hat{a}_{4}^{p}+ \Big(\frac{\iota^{2p+1}}{2p+1}\Big)^{\frac{1}{2}}\hat{a}_{5}^{p}\Big\}\Big\}, \\
&  \oplus_{2}:=5^{p-1}\Big\{M^{p} \ \iota^{p-1} \  \hat{a}_{2}^{p}+ C_{p}\  M^{P}\ \iota^{\frac{p}{2}-1} \ \hat{a}_{3}^{p}\Big\}.
\end{align*}
This implies that 
\begin{align}\label{eq14}
\mathbb{E}\|A(\iota)\|^{p}& \leq \ \frac{\oplus_{2}}{1-\oplus_{1}} \int\limits_{0}^{\iota}(\iota-\upsilon)\mathbb{E}\|A(\upsilon)\|^{p}d\upsilon.
\end{align}
By applying Gronwall's inequality, the above inequality \eqref{eq14} reduces to given Cauchay problem \eqref{eq1}. Also, 
\begin{align*}
\mathbb{E}\|A(\iota)\|^{p} & \ = \ 0\\
 \therefore  \qquad x_{2}(\iota) & \ = \  x_{1}(\iota) \  a.s,
\end{align*}
which means that the mild solution of the give Cauchy problem \eqref{eq1} is unique. Thus, the solution is well defined.
\section{Exponential Stability via Integral Contractors}
In this section, the sufficient criteria of the mild solution for the given Cauchy problem \eqref{eq1} is  investigation by employing the impulsive integral inequality.\\
In order to prove our main result, some additional assumptions are imposed.
\begin{itemize}
\item [$(H_{2})$] For a strongly continuous $ \alpha$-order cosine families $ C_{\alpha}(\iota) $ associated with sine operator $ S_{\alpha}(\iota)$ s.t there exist positive constants $ a_{1} $ and $ a_{2} $ with $ D_{1}, D_{2}>1$ s.t 
$$\sup\limits_{\iota\geq 0} \|C_{\alpha}(\iota) \leq D_{1} e^{-a_{1}\iota}; \  \sup\limits_{\iota\geq 0} \|S_{\alpha}(\iota) \leq D_{2} e^{-a_{2}\iota}.$$
\end{itemize}
\begin{lemma} \label{iie}
	Suppose that for \ $h >0, \eta_{1}, \eta_{2} \in ( 0, h]$,  there exist constants \ $ \xi_{i} > 0 \ (i= 1, 2, 3, 4)$ and a function $ \Psi : [-\kappa, \infty) \rightarrow [0,  \infty)$ s.t
	\begin{align}
	\Psi (\iota)\leq \left\{
	\begin{array}{ll}
	\xi_{1} e^{- \eta_{1} \iota} + \xi_{2} e^{- \eta_{2} \iota}, \quad \iota\in [-r, 0]  \\
	\xi_{1}e^{- \eta_{1} \iota}+\xi_{2} e^{- \eta_{2} \iota}+\xi_{3}\sup\limits_{\theta \in [-r, 0]} \Psi(\upsilon+ \theta)+\xi_{4} \int\limits_{0}^{\iota}e^{-\eta_{1}(\iota-\upsilon)}\sup\limits_{\theta \in [-r, 0]} \Psi(\upsilon+ \theta)d\upsilon \\
	+\xi_{5} \int\limits_{0}^{\iota}e^{-\eta_{2}(\iota-\upsilon)}\sup\limits_{\theta \in [-r, 0]} \Psi(\upsilon+ \theta)d\upsilon+\xi_{6} \int\limits_{0}^{\iota}e^{-\eta_{1}(\iota-\upsilon)}\sup\limits_{\theta \in [-r, 0]} \Psi(\upsilon+ \theta)d\upsilon\\
	+\xi_{7} \int\limits_{0}^{\iota}e^{-\eta_{2}(\iota-\upsilon)}\sup\limits_{\theta \in [-r, 0]} \Psi(\upsilon+ \theta)d\upsilon, \label{z}\quad  t  \geq 0,
	\end{array}
	\right. 
	\end{align}
\end{lemma}
and if
	\begin{align}\label{01}
	\xi_{3}+ \frac{\xi_{4}}{\eta_{1}} -  \frac{\xi_{5}}{\eta_{2}} & < 1
	\end{align}
	then, we have 
	\begin{align}\label{02}
	\Psi (\iota) & \leq N_{\epsilon} e^{-\mu \iota}\  \text{for} \ t \geq -\kappa,
	\end{align} 
	where $\mu \in (0, \eta_{1} \varLambda \ \eta_{2})$ is a positive root of the equation $\xi_{3}e^{-\mu \theta}+\xi_{4} \frac{e^{-\mu \theta}}{\eta_{1} - \mu} + \xi_{5} \frac{e^{-\mu \theta}}{\eta_{2} - \mu} = 1$ and 
	\begin{align*}
	N_{\epsilon} = max  \Big\{\xi_{1} + \xi _{2}, \frac{(\eta_{1} - \mu)}{\xi_{4}  e^{\mu \theta}-\xi_{6}},\frac{(\eta_{2} - \mu)}{\xi_{5}  e^{\mu \theta}-\xi_{7}} \Big\} > 0.  
	\end{align*} 
\begin{theorem}\label{thm2}
	Assume that the assumption $(H_{2})$ is fulfilled, then the given Cauchy problem \eqref{eq1} is exponentially stable in the $p^{th}$ moment sense on $ J, $ provided
	\begin{align*}
	&5^{p-1}\Big[\|\mathcal{A}^{-\gamma}\|^{p}\ \hat{a}_{1}^{p}+\alpha^{p}\ c_{\mu}^{p}\ \hat{a}_{1}^{p}\ \Big[\frac{\iota^{p\alpha \mu}}{p\alpha \mu}\Big]^{p}+D_{2}^{p} a_{2}^{p-1}\ \hat{a}_{2}^{p}
	+C_{p}\  D_{2}^{P}\ \Big( \frac{2a_{2}(p-1)}{p-2}\Big)^{1-\frac{p}{2}} \ \hat{a}_{3}^{p}\\
	&\quad +k(p)\ D_{2}^{p}(\hat{a}_{4}^{\frac{p}{2}}+\hat{a}_{5}^{p})\Big( \frac{2a_{2}(p-1)}{p-2}\Big)^{1-\frac{p}{2}}\Big] <1.
	\end{align*}
	\end{theorem}
\textbf{Proof :} Let $ x(\iota) $ be the mild solution of the given Cauchy problem \eqref{eq1}. Now
\begin{align}
\mathbb{E}\|x(\iota)\|^{p} &\leq \ 5^{p-1}\Big\{\mathbb{E}\Big\|g\Big(\iota, x_{\iota}+ y_{\iota}+\int\limits_{0}^{\iota}S_{\alpha}(\iota-\upsilon)\Gamma_{1}(\upsilon, x_{\upsilon}) y(\upsilon)d\upsilon+\int\limits_{0}^{\iota}S_{\alpha}(\iota-\upsilon)\Gamma_{2}(\upsilon, x_{\upsilon}) y(\upsilon)d\upsilon\nonumber\\
	&\quad+ \int\limits_{0}^{\iota}S_{\alpha}(\iota-\upsilon)\Gamma_{3}(\upsilon, x_{\upsilon}) y(\upsilon)dw(\upsilon)+\int\limits_{0}^{\iota}\int\limits_{\mathcal{Z}}S_{\alpha}(\iota-\upsilon)\Gamma_{4} (\upsilon, x_{\upsilon}, u) y(\upsilon)\lambda du\Big)\nonumber\\
	&\quad-g(\iota, x_{\iota})- \int\limits_{0}^{\iota}S_{\alpha}(\iota-\upsilon)\Gamma_{1}(, x(\upsilon))y(\upsilon)d\upsilon\Big\|^{p}\nonumber
	\end{align}
	\begin{align}\label{eq17}
&\quad+	\mathbb{E}\Big\|\int\limits_{0}^{\iota}\mathcal{A}S_{\alpha}(\iota-\upsilon)g\Big(\upsilon, x_{\upsilon}+y_{\upsilon}+ \int\limits_{0}^{\kappa}S_{\alpha}(\kappa-\upsilon)\Gamma_{1}(\upsilon, x_{\upsilon}) y(\upsilon)d\upsilon+ \int\limits_{0}^{\kappa}S_{\alpha}(\kappa-\upsilon)\Gamma_{2}(\upsilon, x_{\upsilon}) y(\upsilon)d\upsilon\nonumber\\
	&\quad+ \int\limits_{0}^{\kappa}S_{\alpha}(\kappa-\upsilon)\Gamma_{3}(\upsilon, x_{\upsilon}) y(\upsilon)dw(\upsilon)+\int\limits_{0}^{\kappa}\int\limits_{\mathcal{Z}}S_{\alpha}(\kappa-\upsilon)\Gamma_{4} (\upsilon, x_{\upsilon}, u) y(\upsilon)N(d\upsilon, du)\Big)d\upsilon\nonumber\\
	&\quad-\int\limits_{0}^{\iota}\mathcal{A}S_{\alpha}(\iota-\upsilon)g(\upsilon, x_{\upsilon})d\upsilon-\int\limits_{0}^{\iota}S_{\alpha}(\iota-\upsilon)\Gamma_{1}(\upsilon, x_{\upsilon}) y(\upsilon)d\upsilon\Big\|^{p} \nonumber\\
	&\quad+\mathbb{E}\Big\|\int\limits_{0}^{\iota}S_{\alpha}(\iota-\upsilon)f\Big(\upsilon, x_{\upsilon}+y_{\upsilon}+ \int\limits_{0}^{\kappa}S_{\alpha}(\kappa-\upsilon)\Gamma_{1}(\upsilon, x_{\upsilon}) y(\upsilon)d\upsilon+ \int\limits_{0}^{\kappa}S_{\alpha}(\kappa-\upsilon)\Gamma_{2}(\upsilon, x_{\upsilon}) y(\upsilon)d\upsilon\nonumber\\
	&\quad+ \int\limits_{0}^{\kappa}S_{\alpha}(\kappa-\upsilon)\Gamma_{3}(\upsilon, x_{\upsilon}) y(\upsilon)dw(\upsilon)+\int\limits_{0}^{\kappa}\int\limits_{\mathcal{Z}}S_{\alpha}(\kappa-\upsilon)\Gamma_{4} (\upsilon, x_{\upsilon}, u) y(\upsilon)N(d\upsilon, du)\Big)d\upsilon\nonumber\\
			&\quad-\int\limits_{0}^{\iota}S_{\alpha}(\iota-\upsilon)f(\upsilon, x_{\upsilon})d\upsilon-\int\limits_{0}^{\iota}S_{\alpha}(\iota-\upsilon)\Gamma_{2}(\upsilon, x_{\upsilon}) y(\upsilon)d\upsilon\Big\|^{p}\nonumber\\
			&\quad+\mathbb{E}\Big\|\int\limits_{0}^{\iota}S_{\alpha}(\iota-\upsilon)G\Big(\upsilon, x_{\upsilon}+y_{\upsilon}+ \int\limits_{0}^{\kappa}S_{\alpha}(\kappa-\upsilon)\Gamma_{1}(\upsilon, x_{\upsilon}) y(\upsilon)d\upsilon+ \int\limits_{0}^{\kappa}S_{\alpha}(\kappa-\upsilon)\Gamma_{2}(\upsilon, x_{\upsilon}) y(\upsilon)d\upsilon\nonumber\\
				&\quad+ \int\limits_{0}^{\kappa}S_{\alpha}(\kappa-\upsilon)\Gamma_{3}(\upsilon, x_{\upsilon}) y(\upsilon)dw(\upsilon)+\int\limits_{0}^{\kappa}\int\limits_{\mathcal{Z}}S_{\alpha}(\kappa-\upsilon)\Gamma_{4} (\upsilon, x_{\upsilon}, u) y(\upsilon)N(d\upsilon, du)\Big)dw(\upsilon)\nonumber\\
					&\quad-\int\limits_{0}^{\iota}S_{\alpha}(\iota-\upsilon)G(\upsilon, x_{\upsilon})dw(\upsilon)-\int\limits_{0}^{\iota}S_{\alpha}(\iota-\upsilon)\Gamma_{3}(\upsilon, x_{\upsilon}) y(\upsilon)dw(\upsilon)\Big\|^{p}\nonumber\\
				&\quad+\mathbb{E}\Big\|\int\limits_{0}^{\iota}\int\limits_{\mathcal{Z}}S_{\alpha}(\iota-\upsilon)\sigma\Big(\upsilon, x_{\upsilon}+y_{\upsilon}+ \int\limits_{0}^{\kappa}S_{\alpha}(\kappa-\upsilon)\Gamma_{1}(\upsilon, x_{\upsilon}) y(\upsilon)d\upsilon+ \int\limits_{0}^{\kappa}S_{\alpha}(\kappa-\upsilon)\Gamma_{2}(\upsilon, x_{\upsilon}) y(\upsilon)d\upsilon\nonumber\\
						&\quad+ \int\limits_{0}^{\kappa}S_{\alpha}(\kappa-\upsilon)\Gamma_{3}(\upsilon, x_{\upsilon}) y(\upsilon)dw(\upsilon)+\int\limits_{0}^{\kappa}\int\limits_{\mathcal{Z}}S_{\alpha}(\kappa-\upsilon)\Gamma_{4} (\upsilon, x_{\upsilon}, u) y(\upsilon)N(d\upsilon, du)\Big)N(d\upsilon, du)\nonumber\\
						&\quad -\int\limits_{0}^{\iota}\int\limits_{\mathcal{Z}}S_{\alpha}(\iota-\upsilon)\sigma(\upsilon, x_{\upsilon}, u)N(d\upsilon, du)-\int\limits_{0}^{\iota}\int\limits_{\mathcal{Z}}S_{\alpha}(\iota-\upsilon)\Gamma(\upsilon, x_{\upsilon}, u)y(\upsilon)N(d\upsilon, du) \Big\|^{p}\nonumber\Big\}\\
&\leq \ 5^{p-1}\sum\limits_{i=1}^{5} K_{i}.
\end{align}
Here, it is easy to estimate each term of the R.H.S of the above inequality \eqref{eq17}.\\ 
By using Definition \ref{def2.9} (ii), we get the following estimation.
\begin{align*}
K_{1} &= \ \mathbb{E}\Big\|g\Big(\iota, x_{\iota}+ y_{\iota}+\int\limits_{0}^{\iota}S_{\alpha}(\iota-\upsilon)\Gamma_{1}(\upsilon, x_{\upsilon}) y(\upsilon)d\upsilon+\int\limits_{0}^{\iota}S_{\alpha}(\iota-\upsilon)\Gamma_{2}(\upsilon, x_{\upsilon}) y(\upsilon)d\upsilon\nonumber\\
	&\quad+ \int\limits_{0}^{\iota}S_{\alpha}(\iota-\upsilon)\Gamma_{3}(\upsilon, x_{\upsilon}) y(\upsilon)dw(\upsilon)+\int\limits_{0}^{\iota}\int\limits_{\mathcal{Z}}S_{\alpha}(\iota-\upsilon)\Gamma_{4} (\upsilon, x_{\upsilon}, u) y(\upsilon)\lambda du\Big)\nonumber\\
	&\quad-g(\iota, x_{\iota})- \int\limits_{0}^{\iota}S_{\alpha}(\iota-\upsilon)\Gamma_{1}(, x(\upsilon))y(\upsilon)d\upsilon\Big\|^{p}
	\leq \ \|\mathcal{A}^{-\gamma}\|^{p}\ \hat{a}_{1}^{p}\ \mathbb{E}\|x_{\iota}\|^{p}.
	\end{align*}
By using Definition \ref{def2.7} and Definition \ref{def2.9} (ii), we get the following estimation.
\begin{align*}
K_{2} &= \ \mathbb{E}\Big\|\int\limits_{0}^{\iota}\mathcal{A}S_{\alpha}(\iota-\upsilon)g\Big(\upsilon, x_{\upsilon}+y_{\upsilon}+ \int\limits_{0}^{\kappa}S_{\alpha}(\kappa-\upsilon)\Gamma_{1}(\upsilon, x_{\upsilon}) y(\upsilon)d\upsilon+ \int\limits_{0}^{\kappa}S_{\alpha}(\kappa-\upsilon)\Gamma_{2}(\upsilon, x_{\upsilon}) y(\upsilon)d\upsilon\nonumber\\
	&\quad+ \int\limits_{0}^{\kappa}S_{\alpha}(\kappa-\upsilon)\Gamma_{3}(\upsilon, x_{\upsilon}) y(\upsilon)dw(\upsilon)+\int\limits_{0}^{\kappa}\int\limits_{\mathcal{Z}}S_{\alpha}(\kappa-\upsilon)\Gamma_{4} (\upsilon, x_{\upsilon}, u) y(\upsilon)N(d\upsilon, du)\Big)d\upsilon\nonumber\\
	&\quad-\int\limits_{0}^{\iota}\mathcal{A}S_{\alpha}(\iota-\upsilon)g(\upsilon, x_{\upsilon})d\upsilon-\int\limits_{0}^{\iota}S_{\alpha}(\iota-\upsilon)\Gamma_{1}(\upsilon, x_{\upsilon}) y(\upsilon)d\upsilon\Big\|^{p}\\
& \leq \ \alpha^{p}\ c_{\mu}^{p}\ \hat{a}_{1}^{p}\ \Big[\frac{\iota^{p\alpha \mu}}{p\alpha \mu}\Big]^{p}\mathbb{E}\|x_{\iota}\|^{p}.	
\end{align*}
By using assumption $ (H_{2})$ and Definition \ref{def2.9} (i), we get the following estimation. 
\begin{align*}
K_{3} &= \ \mathbb{E}\Big\|\int\limits_{0}^{\iota}S_{\alpha}(\iota-\upsilon)f\Big(\upsilon, x_{\upsilon}+y_{\upsilon}+ \int\limits_{0}^{\kappa}S_{\alpha}(\kappa-\upsilon)\Gamma_{1}(\upsilon, x_{\upsilon}) y(\upsilon)d\upsilon+ \int\limits_{0}^{\kappa}S_{\alpha}(\kappa-\upsilon)\Gamma_{2}(\upsilon, x_{\upsilon}) y(\upsilon)d\upsilon\nonumber\\
&\quad+ \int\limits_{0}^{\kappa}S_{\alpha}(\kappa-\upsilon)\Gamma_{3}(\upsilon, x_{\upsilon}) y(\upsilon)dw(\upsilon)+\int\limits_{0}^{\kappa}\int\limits_{\mathcal{Z}}S_{\alpha}(\kappa-\upsilon)\Gamma_{4} (\upsilon, x_{\upsilon}, u) y(\upsilon)N(d\upsilon, du)\Big)d\upsilon\nonumber\\
&\quad-\int\limits_{0}^{\iota}S_{\alpha}(\iota-\upsilon)f(\upsilon, x_{\upsilon})d\upsilon-\int\limits_{0}^{\iota}S_{\alpha}(\iota-\upsilon)\Gamma_{2}(\upsilon, x_{\upsilon}) y(\upsilon)d\upsilon\Big\|^{p}\\
&\leq \ \mathbb{E}\int\limits_{0}^{\iota}\Big\|S_{\alpha}(\iota-\upsilon)f\Big(\upsilon, x_{\upsilon}+y_{\upsilon}+ \int\limits_{0}^{\kappa}S_{\alpha}(\kappa-\upsilon)\Gamma_{1}(\upsilon, x_{\upsilon}) y(\upsilon)d\upsilon+ \int\limits_{0}^{\kappa}S_{\alpha}(\kappa-\upsilon)\Gamma_{2}(\upsilon, x_{\upsilon}) y(\upsilon)d\upsilon\nonumber\\
&\quad+ \int\limits_{0}^{\kappa}S_{\alpha}(\kappa-\upsilon)\Gamma_{3}(\upsilon, x_{\upsilon}) y(\upsilon)dw(\upsilon)+\int\limits_{0}^{\kappa}\int\limits_{\mathcal{Z}}S_{\alpha}(\kappa-\upsilon)\Gamma_{4} (\upsilon, x_{\upsilon}, u) y(\upsilon)N(d\upsilon, du)\Big)d\upsilon\nonumber
\end{align*}
\begin{align*}
&\quad-\int\limits_{0}^{\iota}S_{\alpha}(\iota-\upsilon)f(\upsilon, x_{\upsilon})d\upsilon-\int\limits_{0}^{\iota}S_{\alpha}(\iota-\upsilon)\Gamma_{2}(\upsilon, x_{\upsilon}) y(\upsilon)d\upsilon\Big\|^{p}\\
&\leq \ D_{2}^{p} a_{2}^{p-1}\ \hat{a}_{2}^{p} \ \int\limits_{0}^{\iota} e^{-a_{2}(\iota-\upsilon)} \mathbb{E}\|x_{\upsilon}\|^{p}d\upsilon.
\end{align*}
We estimate $ K_{4} $ by using  assumption $ (H_{2})$ and Definition \ref{def2.9} (iii), 
\begin{align*}
K_{4} & = \ \mathbb{E}\Big\|\int\limits_{0}^{\iota}S_{\alpha}(\iota-\upsilon)G\Big(\upsilon, x_{\upsilon}+y_{\upsilon}+ \int\limits_{0}^{\kappa}S_{\alpha}(\kappa-\upsilon)\Gamma_{1}(\upsilon, x_{\upsilon}) y(\upsilon)d\upsilon+ \int\limits_{0}^{\kappa}S_{\alpha}(\kappa-\upsilon)\Gamma_{2}(\upsilon, x_{\upsilon}) y(\upsilon)d\upsilon\nonumber\\
&\quad+ \int\limits_{0}^{\kappa}S_{\alpha}(\kappa-\upsilon)\Gamma_{3}(\upsilon, x_{\upsilon}) y(\upsilon)dw(\upsilon)+\int\limits_{0}^{\kappa}\int\limits_{\mathcal{Z}}S_{\alpha}(\kappa-\upsilon)\Gamma_{4} (\upsilon, x_{\upsilon}, u) y(\upsilon)N(d\upsilon, du)\Big)dw(\upsilon)\\
&\quad-\int\limits_{0}^{\iota}S_{\alpha}(\iota-\upsilon)G(\upsilon, x_{\upsilon}) d\upsilon
-\int\limits_{0}^{\iota}S_{\alpha}(\iota-\upsilon)\Gamma_{3}(\upsilon, x_{\upsilon})y(\upsilon) dw(\upsilon)\Big\|^{p}\\
&\ \leq \ C_{p}\ D_{2}^{p} \Big[\int\limits_{0}^{\iota}\Big(\mathbb{E}\Big\|G\Big(\upsilon, x_{\upsilon}+y_{\upsilon}+ \int\limits_{0}^{\kappa}S_{\alpha}(\kappa-\upsilon)\Gamma_{1}(\upsilon, x_{\upsilon}) y(\upsilon)d\upsilon+ \int\limits_{0}^{\kappa}S_{\alpha}(\kappa-\upsilon)\Gamma_{2}(\upsilon, x_{\upsilon}) y(\upsilon)d\upsilon\nonumber\\
&\quad+ \int\limits_{0}^{\kappa}S_{\alpha}(\kappa-\upsilon)\Gamma_{3}(\upsilon, x_{\upsilon}) y(\upsilon)dw(\upsilon)+\int\limits_{0}^{\kappa}\int\limits_{\mathcal{Z}}S_{\alpha}(\kappa-\upsilon)\Gamma_{4} (\upsilon, x_{\upsilon}, u) y(\upsilon)N(d\upsilon, du)\Big)dw(\upsilon)\\
&\quad-\int\limits_{0}^{\iota}S_{\alpha}(\iota-\upsilon)G(\upsilon, x_{\upsilon}) d\upsilon
-\int\limits_{0}^{\iota}S_{\alpha}(\iota-\upsilon)\Gamma_{3}(\upsilon, x_{\upsilon})y(\upsilon) dw(\upsilon)\Big\|^{p}\Big)^{\frac{2}{p}}d\upsilon\Big]^{\frac{p}{2}}\\
&\ \leq \ C_{p}\  D_{2}^{P}\ \Big( \frac{2a_{2}(p-1)}{p-2}\Big)^{1-\frac{p}{2}} \ \hat{a}_{3}^{p}  \int\limits_{0}^{\iota}e^{-a_{2}(\iota-\upsilon)}  \mathbb{E}\|x_{\upsilon}\|^{p}d\upsilon.	
	\end{align*}
By using Lemma \ref{2.8}, Definition \ref{2.6} and Definition \ref{def2.9} (iv),  one can obtain	
\begin{align*}
		K_{5} &\ = \ \int\limits_{0}^{\iota}\int\limits_{\mathcal{Z}}S_{\alpha}(\iota-\upsilon)\sigma\Big(\upsilon, x_{\upsilon}+y_{\upsilon}+ \int\limits_{0}^{\kappa}S_{\alpha}(\kappa-\upsilon)\Gamma_{1}(\upsilon, x_{\upsilon}) y(\upsilon)d\upsilon+ \int\limits_{0}^{\kappa}S_{\alpha}(\kappa-\upsilon)\Gamma_{2}(\upsilon, x_{\upsilon}) y(\upsilon)d\upsilon\nonumber\\
		&\quad+ \int\limits_{0}^{\kappa}S_{\alpha}(\kappa-\upsilon)\Gamma_{3}(\upsilon, x_{\upsilon}) y(\upsilon)dw(\upsilon)+\int\limits_{0}^{\kappa}\int\limits_{\mathcal{Z}}S_{\alpha}(\kappa-\upsilon)\Gamma_{4} (\upsilon, x_{\upsilon}, u) y(\upsilon)N(d\upsilon, du)\Big)N(d\upsilon, du)\\
		&\quad-\int\limits_{0}^{\iota}\int\limits_{\mathcal{Z}}S_{\alpha}(\iota-\upsilon)\sigma(\upsilon, x_{\upsilon}, u) N(du, d\upsilon)- \int\limits_{0}^{\iota}\int\limits_{\mathcal{Z}}S_{\alpha}(\iota-\upsilon)\Gamma_{4}(\upsilon, x_{\upsilon}, u) N(d\upsilon, du)\Big\|^{p}\\
		&\ \leq \ k(p)\ D_{2}^{p}(\hat{a}_{4}^{\frac{p}{2}}+\hat{a}_{5}^{p})\Big( \frac{2a_{2}(p-1)}{p-2}\Big)^{1-\frac{p}{2}}\int\limits_{0}^{\iota}e^{-a_{2}(\iota-\upsilon)}  \mathbb{E}\|x_{\upsilon}\|^{p}d\upsilon.	
	\end{align*}	
These together with equation \eqref{eq17},
\begin{align}\label{eq18}
\mathbb{E}\|x(\iota)\|^{p} & \leq \ 5^{p-1}\Big[\|\mathcal{A}^{-\gamma}\|^{p}\ \hat{a}_{1}^{p}+\alpha^{p}\ c_{\mu}^{p}\ \hat{a}_{1}^{p}\ \Big[\frac{\iota^{p\alpha \mu}}{p\alpha \mu}\Big]^{p}\Big]\mathbb{E}\|x_{\iota}\|^{p}+5^{p-1}\Big\{D_{2}^{p} a_{2}^{p-1}\ \hat{a}_{2}^{p}\nonumber\\
&\quad+C_{p}\  D_{2}^{P}\ \Big( \frac{2a_{2}(p-1)}{p-2}\Big)^{1-\frac{p}{2}} \ \hat{a}_{3}^{p} +k(p)\ D_{2}^{p}(\hat{a}_{4}^{\frac{p}{2}}+\hat{a}_{5}^{p})\Big( \frac{2a_{2}(p-1)}{p-2}\Big)^{1-\frac{p}{2}}\Big\}\nonumber\\
&\times \int\limits_{0}^{\iota} e^{-a_{2}(\iota-\upsilon)} \mathbb{E}\|x_{\upsilon}\|^{p}d\upsilon. 
\end{align}
Then, the above inequality \eqref{eq18} are equivalent to L.H.S of Theorem \ref{thm2}, we have
\begin{align*}
\mathbb{E}\|x(\iota)\|^{p} & \leq \ 5^{p-1}\Big[\|\mathcal{A}^{-\gamma}\|^{p}\ \hat{a}_{1}^{p}+\alpha^{p}\ c_{\mu}^{p}\ \hat{a}_{1}^{p}\ \Big[\frac{\iota^{p\alpha \mu}}{p\alpha \mu}\Big]^{p}\Big]e^{-\eta_{1} \iota}\\
&:= \hat{A}_{1} \ e^{-\eta_{1} \iota},
\end{align*}
where $$\hat{A}_{1}:=5^{p-1}\Big[\|\mathcal{A}^{-\gamma}\|^{p}\ \hat{a}_{1}^{p}+\alpha^{p}\ c_{\mu}^{p}\ \hat{a}_{1}^{p}\ \Big[\frac{\iota^{p\alpha \mu}}{p\alpha \mu}\Big]^{p}\Big].$$
By Lemma \ref{iie} and by equation \eqref{eq18}, we have $\mathbb{E}\|x(\iota)\|^{p} \leq \hat{\oplus}e^{-\eta \iota}, \iota \geq \kappa, \ \eta \in (0, \eta_{1}\Lambda \eta_{2})$, where
$$\hat{\oplus}= \max\Big\{\hat{A}_{1}, \hat{A}_{2}:= 5^{p-1}\Big\{D_{2}^{p} a_{2}^{p-1}\ \hat{a}_{2}^{p}
+C_{p}\  D_{2}^{P}\ \Big( \frac{2a_{2}(p-1)}{p-2}\Big)^{1-\frac{p}{2}} \ \hat{a}_{3}^{p} +k(p)\ D_{2}^{p}(\hat{a}_{4}^{\frac{p}{2}}+\hat{a}_{5}^{p})\Big( \frac{2a_{2}(p-1)}{p-2}\Big)^{1-\frac{p}{2}}\Big\}\Big\}.$$
Here $ \eta $ is a +ve root of the equation $ \hat{A}_{1} \ e^{-\eta_{1} \iota} =1$, $ \hat{A}_{1}= \xi_{3}, $ and $\hat{A}_{2} $ is defined as in above equation.
Hence, the fractional stochastic system \eqref{eq1} is exponentially stable.
\section{Example}
We provide an example in this section for validating the theoretical results. Consider the control problem for higher-order fractional neutral stochastic integro-delay differential system given as follows:
\begin{align}\label{eq21}
 &	^{C}D_{0^{+}}^{\alpha}\Big[Z(x, \iota) + \frac{e^{-8\iota}}{\|(\mathcal{A})^{-\gamma}\|} \ x(\iota)\Big] = \frac{\partial^{2}}{\partial x^{2}} \ Z(x, \iota) + \int\limits_{0}^{\iota} e^{\frac{\iota-s}{2}} \frac{\|Z(x, \iota)\|}{49+ \|Z(x, \iota)\|} d\upsilon \nonumber\\
  &\hspace{0.4cm}+ \frac{e^{\upsilon}}{25+ \|Z(x, \iota)\|} \frac{dw(\iota)}{d\iota}+ \int\limits_{\mathcal{Z}}e^{-6(\iota-\upsilon)}\sigma (\upsilon, x_{\upsilon}), u) N(d\upsilon, du), \ \iota \in J, \   \kappa \in [0, \pi],\nonumber \\
 & z(\iota, 0) = z(\iota, \pi) = 0,  \  \iota \in [0, b], \quad \frac{\partial}{\partial \iota}(z(0, \kappa)) =x_{1}, \  \kappa \in [0,\pi]. 
 \end{align}
Here, $	^{C}D_{0^{+}}^{\frac{5}{3}}$ denotes Caputo partial derivative of order $\alpha =\frac{5}{3} .$  \\ Let $w(\iota)$ denotes the standard Wiener process.\\
Let $\mathcal{A} : \mathbb{H} \rightarrow \mathbb{H}$ defined as
$$\mathcal{A}\omega = \frac{\partial^{2}}{\partial x^{2}}\  \omega(\iota, x), \quad \omega \in \mathcal{D}(\mathcal{A}),$$
where
\begin{align*}
\mathcal{D}((\mathcal{A})^{\frac{5}{3}}) &= \Big\{ \omega \in \mathbb{H}:  \sum\limits_{n = 1}^{\infty}  n^{\frac{5}{3}} \langle f, e_{n}\rangle e_{n}\in \mathbb{H} : \omega(0) = \omega(\pi) = 0 \Big\}.
\end{align*}
Then $ \mathcal{A} $ has spectral representation 
$$ \mathcal{A}\omega =  \sum\limits_{n= 1}^{\infty} - n^{2}(\omega, \omega_{n}),$$
where, $\omega_ {n}(\upsilon)= \sqrt {\frac{2}{\pi}} \  \sin (n \upsilon)$ is the orthonormal set of eigenvalue of $\mathcal{A}$. $\mathcal{A}$ is the infinitesimal generator of a strongly continuous cosine family $\{C(\iota), \iota \in \mathbb{R}\}$ define as 
\begin{align*}
C(\iota)\omega &= \sum\limits_{n=1}^{\infty} cosn\iota <\omega, \omega_{n}> , \ \omega\in \mathbb{H},
\end{align*}
 and the associated sine family is given by 
 \begin{align*}
 S(\iota)\omega &= \sum\limits_{n=1}^{\infty}- \frac{1}{n} sin n\iota <\omega, \omega_{n}> , \ \omega\in \mathbb{H}.
 \end{align*}
  Hence, from the subordinate principle, it follows that $\mathcal{A}$ is the exponentially bounded fractional cosine family $C(\iota)$ such that $C(0) = I.$ The nonlinear functions $f: J \times \mathcal{B}_{r} \rightarrow \mathbb{H}, g: J \times \mathcal{B}_{r} \rightarrow \mathbb{H}, G: J \times \mathcal{B}_{r} \rightarrow L_{Q}^{0}(\mathbb{K},\mathbb{H})$, and $\sigma:J \times \mathcal{B}_{r} \times \mathcal{Z}\rightarrow \mathbb{H}, $  described by 
\begin{align*}
g(\iota, x_{\iota})& = \ \frac{e^{-2\iota}}{\|\mathcal{A}^{-\gamma}\|},\\
 f(\iota, x_{\iota}) &= \ \frac{\frac{e^{\iota-s}}{2}}{49},\\
G(\iota, x_{\iota})\frac{dw(\iota)}{d\iota} &= \  \frac{e^{\upsilon}}{25} \frac{dw(\iota)}{d\iota},\\
\int\limits_{\mathcal{Z}}\sigma(\iota, x_{\iota}, u) N(d\iota, du)&= \int\limits_{\mathcal{Z}} e^{-6(\iota-\upsilon)}N(d\iota, du). 
\end{align*}
The non-linear functions are $ f, g, G,$  $\sigma $ have  regular integral contractors $\Gamma_{1}= \Gamma_{2}=\Gamma_{3}=\Gamma_{4} = 0.$

\noindent
Hence, the system \eqref{eq1} has a unique mild solution. By choosing the particular values in the given parameters are, $M=0.002, \ p= 2,\ \iota = \frac{1}{2}, C_{p}= 1, \hat{a}_{1}= \hat{a}_{2}= \hat{a}_{3}=\hat{a}_{4}= \frac{1}{3}$, and $\hat{a}_{5}= 0.02, \  k(p)=1.$\\

From the definition of $(\mathcal{A})^{-\gamma}$ [see \cite{r13}], it is easy to infer that
$$\|(\mathcal{A})^{-\gamma}\| \leq \frac{1}{\pi^{\frac{3}{2}}}.$$
Accordingly, all the hypotheses of Theorem \ref{PIDthm1} are satisfied.
\begin{align*}
&5^{p-1} \Big\{\|\mathcal{A}^{-\gamma}\|^{p} \ \hat{a}_{1}^{p}+\alpha^{p}\ c_{\mu}^{p}\ \hat{a}_{1}^{p}\ \Big[\frac{\iota^{p\alpha \mu}}{p\alpha \mu}\Big]^{p}+k(p)\ M^{p}\Big\{ \Big(\frac{\iota^{3}}{3}\Big)^{\frac{p}{2}}\hat{a}_{4}^{p}+ \Big(\frac{\iota^{2p+1}}{2p+1}\Big)^{\frac{1}{2}}\hat{a}_{5}^{p}\Big\} \nonumber\\
	&\quad+M^{p} \ \iota^{p-1} \  \hat{a}_{2}^{p}+ C_{p}\  M^{P}\ \iota^{\frac{p}{2}-1} \ \hat{a}_{3}^{p}\Big\}\\
&:=5^{p-1}\Big\{ e^{-1} \frac{1}{\pi^{\frac{3}{2}}} + \Big(\frac{5}{3}\Big)^{2} \Big(\frac{1}{3}\Big)^{2}\Big[\frac{\frac{1}{2}2 \times 0.5} {2 \times 0.5}\Big]^{2}+1 \times 0.002 \times 1 \times \frac{e^{0.5}}{25}\\
&\quad+1 \times 0.002\Big(\frac{(0.5)^{3}}{3}\times e^{-6\times 0.5}+\Big(\frac{(0.5)^{5}}{5}\Big)^{\frac{1}{2}} e^{-6\times 0.5}\Big) \Big\}\\
	 &:=0.007436< \ 1.
\end{align*}
Moreover, by Theorem \ref{thm2}, we may deduce that if
\begin{align*}
&5^{p-1}\Big[\|\mathcal{A}^{-\gamma}\|^{p}\ \hat{a}_{1}^{p}+\alpha^{p}\ c_{\mu}^{p}\ \hat{a}_{1}^{p}\ \Big[\frac{\iota^{p\alpha \mu}}{p\alpha \mu}\Big]^{p}+D_{2}^{p} a_{2}^{p-1}\ \hat{a}_{2}^{p}
	+C_{p}\  D_{2}^{P}\ \Big( \frac{2a_{2}(p-1)}{p-2}\Big)^{1-\frac{p}{2}} \ \hat{a}_{3}^{p}\\
	&\quad +k(p)\ D_{2}^{p}(\hat{a}_{4}^{\frac{p}{2}}+\hat{a}_{5}^{p})\Big( \frac{2a_{2}(p-1)}{p-2}\Big)^{1-\frac{p}{2}}\Big] \\
	&:=5^{p-1}\Big\{ e^{-1} \frac{1}{\pi^{\frac{3}{2}}} + \Big(\frac{5}{3}\Big)^{2} \Big(\frac{1}{3}\Big)^{2}\Big[\frac{\frac{1}{2}2 \times 0.5} {2 \times 0.5}\Big]^{2}+1 \times 0.002 \times 1 \times 0.33\\
	&\quad+1 \times 0.002\Big(\frac{(0.5)^{3}}{3}\times e^{-6\times 0.5}+\Big(e^{-6 \times 0.5}(0.33) \Big)\Big) \Big\}\\
	& := 0.00874<1.
\end{align*}
Then the mild solution of the system \eqref{eq1} is  exponentially stable.
\section{Conclusion}
A new exponential stability model is presented with existence and uniqueness for the higher-order fractional neutral delay differential system using the new technique of integral contractors with the regularity.
To demonstrate the results, stochastic analysis approach and the concept of bounded integral contractors was combined with the sequencing technique. Furthermore, exponential stability results for fractional stochastic equations have been established by using impulsive integral inequality technique. The numerical example, it helps to establish the results numerically with simulation and one can give an application in the numerical part of exponential stability using this result. In future, authors have planned the same phenomena to study fractional stochastic partial differential equation models. Caputo fractional derivative can also be replaced by Hilfer fractional derivative as of future study. \\

{\bf{Conflict of Interest:}} Authors declare that they do not have any conflict of interest.

\end{document}